\documentclass[11pt,a4paper,reqno]{amsart}%
\usepackage{amsthm,amsmath,amsfonts,amssymb,xcolor,amsxtra,appendix,bookmark,dsfont,bbm,mathrsfs,enumitem,float}
\usepackage{makecell}
\usepackage[latin1]{inputenc}
\usepackage{color} 
\usepackage{amssymb}\usepackage{graphicx}
\usepackage{tikz}
\usepackage{hyperref}
\theoremstyle{plain}
\newtheorem{theorem}{Theorem}%[section]
\newtheorem{lemma}[theorem]{Lemma}

\newtheorem{proposition}[theorem]{Proposition}

\newtheorem{definition}{Definition}

\theoremstyle{remark}
\newtheorem{remark}[theorem]{Remark}

\newcommand{\bbd}{\mathbbm d}

\newcommand{\R}{\mathbb R}
\newcommand{\N}{\mathbb N}
\newcommand{\C}{\mathbb C}
\newcommand{\Z}{\mathbb Z}
\newcommand{\T}{\mathbb T}

\newcommand{\E}{\mathbb E}

\newcommand{\sS}{\mathscr S}

\newcommand{\sN}{\mathscr N}

\DeclareMathOperator{\Hess}{Hess}
\DeclareMathOperator{\diag}{diag}

\newcommand\numberthis{\addtocounter{equation}{1}\tag{\theequation}}

\allowdisplaybreaks

\title[Internal control of the transition kernel for stochastic lattice dynamics]{Internal control of the transition kernel for stochastic lattice dynamics}
\author[A. Hannani]{Amirali Hannani
}
\address{Institute for Theoretical Physics, KU Leuven, 3001 Leuven, Belgium}
\email{amirali.hannani@kuleuven.be
} 
\author[M. N. Phung]{Minh Nhat Phung}
\address{Department of Mathematics, Texas A\&M University, College Station, TX 77843, USA}
\email{pmnt1114@tamu.edu}

\author[M.-B. Tran]{Minh-Binh Tran}
\address{Department of Mathematics, Texas A\&M University, College Station, TX 77843, USA}
\email{minhbinh@tamu.edu} 
\thanks{M.-B. T is  funded in part by  a   Humboldt Fellowship,   NSF CAREER DMS-2303146, and NSF Grants DMS-2204795, DMS-2305523, DMS-2306379.}

\author[E. Tr\'elat]{Emmanuel Tr\'elat}
\address{Sorbonne Universit\'e, CNRS, Universit\' e de Paris, Inria, Laboratoire Jacques-Louis Lions (LJLL),F-75005 Paris, France}
\email{emmanuel.trelat@sorbonne-universite.fr} 
\begin{document}

	\date{\today}

	\begin{abstract} In   \cite{hannani2024controlling}, we have designed impulsive and feedback controls for   harmonic chains with a point thermostat. In this work, we study the internal control for stochastic lattice dynamics, with the goal of controlling the transition kernel of the kinetic equation in the limit. A major novelty of the work is the introduction of a new geometric combinatorial argument, used to establish paths for the controls.
	\end{abstract}
	
	\maketitle

	\tableofcontents

\section{Introduction}
During the last few years, energy transport in anharmonic chains, including the famous Fermi-Pasta-Ulam (FPU) chain with pinning potential and a quartic nonlinearity, has become an important topic of research (see \cite{basile2010energy,komorowski2023heat,komorowski2022asymptotic,KORH2020}). The Hamiltonian of the FPU-chain is defined by
\begin{equation}\label{FPU}
	H(\beta,\alpha):=\sum_{n=-N}^{N}\left( \frac{1}{2M}\alpha_n^2+\frac{1}{2}\omega_0^2\beta_n^2 \right)+\sum_{n=-N}^{N}\left( \frac{1}{2}A(\beta_{n+1}-\beta_n)^2+B(\beta_{n+1}-\beta_n)^4 \right)
\end{equation}
where $\beta_n$ is the displacement of the particle $n$ from its equilibrium position, $\alpha_n$ is the canonically conjugate momentum and $M$ is the mass of each particle. The pinning potential represented by the term $\omega_0^2\beta_n^2$ contributes to confining the particle $n$. In \eqref{FPU}, the second sum depending only on the displacement differences consists of a harmonic nearest neighbor term $\frac{1}{2}A(\beta_{n+1}-\beta_n)^2\ (A>0)$ and an anharmonic term $B(\beta_{n+1}-\beta_n)^4\ (B\ge0)$. The mathematical study of energy transport in the FPU-chain is challenging.

The authors of \cite{basile2010energy} replaced the nonlinearity in \eqref{FPU} by a stochastic exchange of momentum between nearest neighboring sites chosen such that the exchange conserves the local energy. In another word, they set $B=0$ and add a small perturbation to the system.

When $B=0$, \eqref{FPU} corresponds to a discrete wave equation. The energy transport of the wave equation has been studied extensively using Wigner distributions (see \cite{basile2010energy,KORH2020,LS2007,RPK96}). When the wave propagates without ``noise'', the Wigner distribution is governed by the linear transport equation
\begin{equation}\label{nonoise}
	\frac{d}{dt}W(x,k,t)+\frac{1}{2\pi}\omega'(k)\frac{\partial}{\partial x}W(x,k,t)=0
\end{equation}
where $x\in\R$ is the position along the chain, the wave number $k$ belongs to $\T$, the unit torus identified with $\left[ -\frac{1}{2},\frac{1}{2} \right]$ with periodic endpoints. The variable $t$ is the time. The dispersion relation inferred from \eqref{FPU} is 
\begin{equation}\label{omegaex}
	\omega(k)^2=\omega_0^2+2A(1-\cos(2\pi k)).
\end{equation}
Indeed, in \eqref{FPU}, the terms involving $\beta$ can be rewritten as
\begin{equation}
\label{FPUbeta}
\frac{1}{2}\sum_{n\in\Z}\sum_{m\in\Z}\beta_m\beta_n(\omega_0^2\delta_{m-n}+2A\delta_{m-n}-A\delta_{m-n-1}-A\delta_{m-n+1}),
\end{equation}
where $\delta_n$ is the Kronecker delta function.
The Fourier transform of \eqref{FPUbeta} is
\begin{equation*}
	\frac{1}{2}\hat{\beta}(k)^2(\omega_0^2+2A(1-\cos(2\pi k))).
\end{equation*}
The coefficient of $\hat{\beta}^2$ gives the dispersion relation \eqref{omegaex}.
Equation \eqref{nonoise} describes the energy transport because the local energy density time $t$ is defined either by
\begin{equation*}
	E(k,t):=\int_{\R}W(x,k,t)dx \quad\text{or by}\quad\bar{E}(x,t):=\int_{\T}W(x,k,t)dk.
\end{equation*}
When we take into consideration a ``noise'' in the discrete wave equation, we get
\begin{equation}\label{sysmodel}
	\begin{split}
		d\beta_n(t)=& \alpha_n(t)dt,\\
		d\alpha_n(t)=& -\sigma\star \beta_n(t)dt+P^{\varepsilon}_{N,n}(\alpha,dt),
	\end{split}
\end{equation}
where $\sigma$ is a coupling between two particles $n,n'$ depending only on their distance $|n-n'|$, $\sigma\star\beta_n=\sum_{n'\in\Z}\sigma_{n-n'}\beta_n'$ is the discrete convolution, and $P^{\varepsilon}_{N,n}(\alpha,dt)$ is the perturbation corresponding to the stochastic exchange of momentum between neighboring particles within distance $N$ of the particle $n$ conserving the local energy. With the perturbation, \eqref{nonoise} becomes
\begin{equation}\label{introker}
	\frac{d}{dt}W(x,k,t)+\frac{1}{2\pi}\omega'(k)\frac{\partial}{\partial x}W(x,k,t)=\int_{\T}\mathcal{K}(k,k')(W(x,k',t)-W(x,k,t))dk,
\end{equation}
where $\mathcal{K}$ is a \emph{transition kernel}. We refer to \cite{borcea2016derivation,pomeau2019statistical} for  related situations.

Control problems in kinetic theory have become an important topic  in the recent years  with significant advances(see \cite{carrillo2022controlling,piccoli2015control}). Relating to this progress, we designed impulsive and feedback controls for the harmonic chains with a point thermostat in our previous work \cite{hannani2024controlling}. In the current work, the proofs are built base on a novel geometric combinatorial argument.

Based on    \cite{basile2010energy},    we address the internal control problem for the harmonic chain with perturbation, with the purpose of controlling the  transition kernel. For a given $\hat{K}:\T\times\T\to\R$, we want to find $P^\varepsilon_{N,n}$ that leads to \eqref{introker} with $\mathcal{K}\approx\hat{K}$. For this purpose, we design stochastic exchanges of momentum between neighboring particles of distance $N\in\mathbb{N}$. This is significantly harder than the  model considered in \cite{basile2010energy} where $N$ is equal to $1$. Considering a larger $N$ increases the number of free parameters and thus the number of controls, which is a significant improvement in terms of control theory.  Furthermore,  the $N$-neighboring problem is then extended to the multi-dimensional case. To address the  $N$-neighboring problem in the multi-dimensional case, we introduce the new concept of simple index and dual indices for the exchanges of momentum between neighboring particles. 

Let us make an additional remark on the transition kernel. Since the kernel satisfies
\begin{equation*}
	\int_{\T}\mathcal{K}(k,k')dk'\sim |k|^2\quad\text{for }|k|\ll 1,
\end{equation*}
we consider target kernels of the form 
\begin{equation*}
	\sin^2(\pi k)\sin^2(\pi k')\tilde{v}(\cos(\pi k),\cos(\pi k')),
\end{equation*}
where $\tilde{v}$ is a polynomial with two variables. In our work, we compute and characterize the set $V$ of possible achievable polynomials $\tilde{v}$. For each $\tilde{v}\in V$, we design a control achieving the resulting kernel.

To be more specific about the nature of our control,  we set up $2N$ controlled parameters $M_N$. Among them we can freely choose values for up to $N-\lfloor N/2\rfloor$ parameters. Then, there exists a perturbation and a kernel corresponding to the chosen parameters. The resulting kernel $\hat{K}_{M_N}$ yields the linear transport equation \eqref{introker} for the Wigner distribution. Hence, by choosing appropriate parameters, we design an explicit control for the energy transport problem.

We also extend our results to the multidimensional torus $\mathbb{T}^\bbd$ . On the space of dimension $\bbd\ge 2$, we introduce controlled parameters $M_{\gamma}$ for simple index and $M^{a,b}_{\gamma}$ for dual indices. In the multi-dimensional case, we actually have more free parameters. In the simple index case, we can have up to $2^{\bbd-1}\bbd^N(N-\lfloor N/2\rfloor)$ free parameters. For dual indices and $N,\bbd>1$, we have up to $2^{\bbd-1}\bbd^{N+1}(\bbd-1)(N-1)$ free parameters.

Finally, to clarify our computational steps, we give an example for small distance. In the appendix, we perform detailed computations for $K_{M_N}$ with $N=3$ and then show explicitly the achievable target kernels.

\section{The Setting}

\subsection{Preliminaries}
\label{Prelim}
We normalize the mass $M=1$ and consider the harmonic case $B=0$. The Hamiltonian \eqref{FPU} is then
\begin{equation*}
	H(\beta,\alpha)=\frac{1}{2}\sum_{n\in\Z}\alpha_n^2+\frac{1}{2}\sum_{n,n'}\sigma_{n-n'}\beta_n\beta_{n'}.
\end{equation*}
For a function $f \in \ell^2(\Z)$, the Fourier transform is
\begin{equation*}
	\hat{f}(k):=\sum_{n\in\Z}e^{-2\pi ikn}f_n,\quad\forall k\in\T.
\end{equation*}
Also, the inverse Fourier transform of a function $f\in L^2(\T)$ is
\begin{equation*}
	\tilde{f}_n:=\int_{\T}e^{2\pi ikn}f(k)dk.
\end{equation*}
We use these notations for the Fourier transform of the wave functions $\hat{\alpha},\hat{\beta},\hat{\phi}$, the coupling $\hat{\sigma}$ and the inverse Fourier transform of dispersion relation $\tilde{\omega}$. The wave $\phi$ is given by
\begin{equation*}
	\phi_n(t):=\frac{1}{\sqrt{2}}\left( \tilde{\omega}\star\beta_n(t)+i\alpha_n(t) \right)
\end{equation*}
or by its Fourier transform
\begin{equation}
	\hat{\phi}(k,t):=\frac{1}{\sqrt{2}}\left( \omega \hat{\beta}(k,t)+i\hat{\alpha}(k,t) \right).
	\label{phidef}
\end{equation}
The dispersion relation is defined by
\begin{equation*}
	\omega(k):=\sqrt{\hat{\sigma}(k)}.
\end{equation*}
We will specify the assumptions on the coupling $\sigma$ later. An example for the dispersion relation is shown in \eqref{omegaex}.

\subsubsection{The setting of \cite{basile2010energy}}
Before going into more details, let us briefly report on the contruction of the perturbation $P_{N,n}^{\varepsilon}$ in \eqref{sysmodel} that was done in \cite{basile2010energy}. The authors considered the $1$-distance perturbation with friction. In this case, the vector field that conserves local momentum and energy is
\begin{equation*}
	\lambda_n:=(\alpha_n-\alpha_{n+1})\frac{\partial}{\partial\alpha_{n-1}}+(\alpha_{n+1}-\alpha_{n-1})\frac{\partial}{\partial\alpha_{n}}+(\alpha_{n-1}-\alpha_n)\frac{\partial}{\partial\alpha_{n+1}}
\end{equation*}
and the corresponding generator of the system \eqref{sysmodel} is
\begin{equation*}
	\sum_{n\in\Z}\alpha_n\frac{\partial}{\partial \beta_n}-\sum_{n,n'\in\Z}\sigma_{n-n'}\beta_{n'}\frac{\partial}{\partial \alpha_n}+\frac{\varepsilon c}{6}\sum_{n \in\Z}(\lambda_n)^2,
\end{equation*}
where $c$ is the friction.
In \eqref{sysmodel}, the perturbation used in \cite{basile2010energy} is
\begin{equation*}
	\begin{split}
		\frac{\varepsilon c}{6}&\Delta(4\alpha_n+\alpha_{n-1}+\alpha_{n+1})dt+\sqrt{\frac{\varepsilon c}{3}}\sum_{d=-1,0,1}(\lambda_{n+d}\alpha_n)Y_{n+d}(dt)\\
		&= \frac{\varepsilon c}{3}(\alpha_{n+2}+\alpha_{n-2}+2\alpha_{n+1}+2\alpha_{n-1}-6\alpha_{n})+\sqrt{\frac{\varepsilon c}{3}}\sum_{d=-1,0,1}(\lambda_{n+d}\alpha_n)Y_{n+d}(dt).
	\end{split}
\end{equation*}
Here, the $Y_n$'s are independent Wiener processes; we will introduce them more formally in Sections \ref{1Dno} and \ref{mulDno}.

At the limit $\varepsilon\to0$, we get the Boltmann transport equation \eqref{introker} with the transition kernel
\begin{equation*}
	\frac{4 c}{3}\left( 2\sin^2(2\pi k)\sin^2(\pi k')+2\sin^2(\pi k)\sin^2(2\pi k')-\sin^2(2\pi k)\sin^2(2\pi k') \right).
\end{equation*}
{\it The purpose of our work is to design a control on this transition kernel by modifying the vector field $\lambda_{n}$ and the perturbation $P_{N,n}^{\varepsilon}$.}

\subsubsection{The 1-dimensional case}\label{1Dno}
The system that governs the wave is \eqref{sysmodel}. The underlying Hamiltonian operator is
\begin{equation}
	O^H:=\sum_{n\in\Z}\alpha_n\frac{\partial}{\partial \beta_n}-\sum_{n,n'\in\Z}\sigma_{n-n'}\beta_{n'}\frac{\partial}{\partial \alpha_n}.
	\label{O^H}
\end{equation}
We apply $O^H$ on \eqref{phidef} and get
\begin{equation}
	O^H\hat{\phi}(k,t)=-i\omega(k)\hat{\phi}(k,t).
	\label{O^Hiden}
\end{equation}

Now, we state the assumptions on the coupling $\sigma$.
\begin{enumerate}[label=($\sigma$\arabic*)]
	\item \label{sig1}There is $n\ne 0$ such that $\sigma_{n}\ne0$.
	\item $\sigma$ is an even function, i.e. $\sigma_{n}=\sigma_{-n}, \forall n\in\Z$.
	\item $\sigma$ decays exponentially, i.e. there exist $C_1,C_2>0$
		\begin{equation*}
			|\sigma_n|\le C_1e^{-C_2|n|},\qquad\forall n\in\Z.
		\end{equation*}
	\item \label{sig4}We have $\hat{\sigma}(k)>0$ for all $k\in\T$ (pinning) or $\hat{\sigma}(k)>0$ for $k\in\T\setminus\{0\}$, $\hat{\sigma}(0)=0$, $\hat{\sigma}''(0)>0$ (no-pinning).
\end{enumerate}

We denote $\{Y_n(t),t\ge0\}_{n\in\Z}$ be a sequence of independent Wiener processes on a probability space with proper filtration $(\Omega,\mathfrak{F}_t,\mathbb{P})$ and $\langle \cdot\rangle_{\mu_{\varepsilon}}$ is the expectation with respect to the probability $\mu_{\varepsilon}$ over the initial state of the wave. We denote by $\E_{\varepsilon}$ the expectation when we need to consider both the initial state $\mu_\varepsilon$ and the Wiener processes with probability $\mathbb{P}$.

The wave function $\phi$ is defined in \eqref{phidef}, where $\alpha,\beta$ are governed by a system with Wiener process. We will specify the system later in \eqref{system} and work with the initial state first.
The initial state probability $\mu_{\varepsilon}$ is assumed to satisfy the following assumptions
\begin{enumerate}[label=($\mu$\arabic*)]
	\item \label{mean0}$\langle \phi_n(0)\rangle_{\mu_\varepsilon}=0, \forall n\in\Z$.
	\item $\langle \phi_n(0)\phi_{n'}(0)\rangle_{\mu_\varepsilon}=0, \forall n,n'\in\Z$.
	\item\label{boundenergy} There is $C_3>0$ such that $\varepsilon\langle\|\phi(0)\|_{\ell^2}^{2}\rangle_{\mu_{\varepsilon}}<C_3, \forall \varepsilon>0$.
	\item\label{nopin} In the no-pinning case,
		\begin{equation*}
			\lim_{R\to 0}\limsup_{\varepsilon\to0}\frac{\varepsilon}{2}\int_{|k|<R}\langle|\hat{\phi}(k,0)|^2\rangle_{\mu_{\varepsilon}}dk=0.
		\end{equation*}
\end{enumerate}

The Wigner distribution is denoted by $W^{\varepsilon}$ and is defined by
\begin{equation*}
	\begin{split}
		\langle S,W^{\varepsilon}\rangle:=& \frac{\varepsilon}{2}\sum_{n,n'\in\Z}\E_{\varepsilon}\left[ \phi_n(t/\varepsilon)\phi_{n'}^{*}(t/\varepsilon) \right]\int_{\T}e^{2\pi ik(n'-n)}S^{*}\left( \frac{\varepsilon(n+n')}{2},k \right)dk\\
		=& \frac{\varepsilon}{2}\sum_{n,n'\in\Z}\E_{\varepsilon}\left[ \phi_{n}(t/\varepsilon)\phi_{n'}^{*}(t/\varepsilon) \right]\tilde{S}^{*}\left( \frac{\varepsilon(n+n')}{2},n-n' \right),
	\end{split}
\end{equation*}
where $S$ is a test function in the Schwartz space $\sS(\R\times\T)$ and $\tilde{S}$ is the inverse Fourier transform
\begin{equation*}\tilde{S}(x,n)=\int_{\T}e^{2\pi ikn}S(x,k)dk.\end{equation*}
For the initial state, we also assume that
\begin{enumerate}[label=(W\arabic*)]
	\item \label{mu0}$W^\varepsilon(dx,dk,0)$ converges (at least, a subsequence of it) to a non-negative measure $\mu_0(dx,dk)$.
\end{enumerate}
We also define the energy density $E^{\varepsilon}(t)$ as a measure on $\T$
\begin{equation*}
	\langle S,E^{\varepsilon}(t)\rangle:=\langle S,W^{\varepsilon}(t)\rangle= \frac{\varepsilon}{2}\int_{\T}\E_{\varepsilon}[|\hat{\phi}(k,t/\varepsilon)|^2]S(k)dk
\end{equation*}
where $S(k)$ is a bounded real-valued function on $\T$. According to Assumption \ref{mu0}, we also have that $E^\varepsilon(dk,0)$ converges to a measure $\nu_0(dk)$.

By \ref{boundenergy} and the conservation of energy, we can show that the Wigner distribution is well-defined. In fact, there is $C_4>0$ such that
\begin{equation*}
	\begin{split}
		|\langle S,W^{\varepsilon}\rangle|&\le \frac{\varepsilon}{2}\E_{\varepsilon}\left[ \|\phi(t/\varepsilon)\|_{\ell^2}^2 \right]\sum_{n}\sup_{x\in\R}\left|\tilde{S}^{*}\left( x,n \right)\right|\le C_4.
	\end{split}
\end{equation*}
On the space $\sN$ defined as the completion of the Schwartz space for the norm
\begin{equation*}
	\|S\|_{\sN}:=\sum_{n\in\Z}\sup_{x\in\R}|\tilde{S}(x,n)|,
\end{equation*}
$W^{\varepsilon}$ is a continuous linear operator, and
\begin{equation*}
	\|W^{\varepsilon}\|_{\sN'}<\infty.
\end{equation*}
This also implies that the family $(W^{\varepsilon})_{0<\varepsilon<1}$ is sequentially weak-* compact. We can take a subsequence, if it is necessary, to get Assumption \ref{mu0}.

\subsubsection{The multi-dimensional case}\label{mulDno}
In the multi-dimensional case, the particles now are placed on the lattice $\Z^{\bbd}$ and $\alpha_{n},\beta_{n}\in\R^{\bbd}$ for each $n\in\Z^{\bbd}$.
The conditions for $\sigma$ and $\mu_{\varepsilon}$ in multi-dimensional are the same as the 1-dimensional case with $\Z$ changed to $\Z^{\bbd}$, $\T$ changed to $\T^{\bbd}$ and $\hat{\sigma}''(0)>0$ changed to $\Hess(\hat{\alpha})(0)$ which is invertible, as follows:
\begin{enumerate}[label=($\sigma$\arabic*')]
	\item \label{sig1'}There is $n \in \Z^\bbd\setminus\{0\}$ such that $\sigma_{n}\ne0$.
	\item $\sigma$ is an even function, i.e. $\sigma_{n}=\sigma_{-n}$ for $n\in\Z^{\bbd}$.
	\item $\sigma$ decays exponentially, i.e. there exist $C_5,C_6$ such that
		\begin{equation*}
			|\sigma_n|\le C_5e^{-C_6|n|},\qquad\forall n\in\Z^{\bbd}.
		\end{equation*}
	\item \label{sig4'}We have $\hat{\sigma}(k)>0$ for all $k\in\T^{\bbd}$ (pinning) or $\hat{\sigma}(k)>0$ for $k\in\T^{\bbd}\setminus\{0\}$, $\hat{\sigma}(0)=0$, $\Hess\hat{\sigma}(0)$ is invertible (no pinning).
\end{enumerate}
In addition, we have the wave functions $\alpha_n,\beta_n\in\R^\bbd$ and $\phi_n\in\C^\bbd$. We write the $a$-th component of those vectors by $[\alpha_n]_a,[\beta_n]_a$ and $[\phi_n]_a$.

We will use the notations $\{Y_{n}^a(t),t\ge0\}_{1\le a\le \bbd,n \in \Z^\bbd}$ (the simple index case in Section \ref{consimple}) and $\{Y^{a,b}_{n,\gamma}(t),t\ge0\}_{1\le a,b\le\bbd,n \in \Z^\bbd,\gamma}$ (the dual indices case in Section \ref{condual}) for sequences of independent Wiener processes. The notation $\gamma$ is a path, this is instrumental for the controls in the multi-dimensional case in Section \ref{Highd}.  Other notations, $\mathbb{P},\mu_\varepsilon,\E_\varepsilon$, are the same as the 1-dimensional case.

We will also specify the system governing the wave in \eqref{systemsimple} or \eqref{systemij}. The assumptions on the initial state are:
\begin{enumerate}[label=($\mu$\arabic*')]
	\item \label{mean0'}$\langle \phi_n(0)\rangle_{\mu_\varepsilon}=0, \forall n\in\Z^{\bbd}$.
	\item $\langle [\phi_n(0)]_a[\phi_{n'}(0)]_b\rangle_{\mu_\varepsilon}=0, \forall n,n'\in\Z^{\bbd},1\le a,b\le \bbd$.
	\item There is $C_7>0$ such that $\varepsilon^{\bbd}\langle\|\phi(0)\|_{\ell^2}^{2}\rangle_{\mu_{\varepsilon}}<C_7, \forall \varepsilon>0$.
	\item \label{nopin'}If the coupling is the no pinning case, then
		\begin{equation*}
			\lim_{R\to 0}\limsup_{\varepsilon\to0}\left(\frac{\varepsilon}{2}\right)^{\bbd}\int_{|k|<R}\langle\|\hat{\phi}(k,0)\|_{L^2(\T^\bbd)}^2\rangle_{\mu_{\varepsilon}}dk=0.
		\end{equation*}
\end{enumerate}

The general Wigner distribution is defined by
\begin{align*}
	\langle J,W^{\varepsilon}(t)\rangle:=& \left( \frac{\varepsilon}{2} \right)^{\bbd}\int_{\R^\bbd\times\T^\bbd}\E_{\varepsilon}[\hat{\phi}(k-\varepsilon\xi/2,t/\varepsilon)\cdot S^{*}(\xi,k)\hat{\phi}(k+\varepsilon\xi/2,t/\varepsilon)]d\xi dk\\
	=& \left( \frac{\varepsilon}{2} \right)^{\bbd}\sum_{n,n'\in\Z^{\bbd}}\sum_{1\le a,b\le \bbd}\E_\varepsilon[[\phi_n]_a[\phi_{n'}^*]_b]\int_{\T^\bbd}e^{i2\pi k\cdot(n'-n)}S^{*}_{b,a}\left(\frac{\varepsilon(n+n')}{2},k\right)dk\\
	=& \left( \frac{\varepsilon}{2} \right)^{\bbd}\sum_{n,n'\in\Z^{\bbd}}\sum_{1\le a,b\le \bbd}\E_\varepsilon[[\phi_n]_a[\phi_{n'}^*]_b]\tilde{S}^{*}_{b,a}\left(\frac{\varepsilon(n+n')}{2},n-n'\right),
\end{align*}
where $S=(S_{a,b})_{1\le a,b\le\bbd}\in\sS(\R^\bbd\times\T^\bbd,\mathbb M_\bbd)$. As in the 1-dimensional case, $W^\varepsilon$ is well-defined and sequentially weak-* compact. The Wigner distribution can be understood as the matrix of distribution, we write it as (see \cite{hannani2022wave,staffilani2021wave} for related ideas)
\begin{equation*}
	W^{\varepsilon}_{b,a}=\left( \frac{\varepsilon}{2} \right)^{\bbd}\int_{\R^\bbd}e^{i2\pi x\cdot\xi}\E_{\varepsilon}[\hat{\phi}^{*}]_b(k-\varepsilon\xi/2,t/\varepsilon)[\hat{\phi}]_a(k+\varepsilon\xi/2,t/\varepsilon)d\xi.
\end{equation*}
We also have a similar assumption to \ref{mu0} for the multi-dimensional case:
\begin{enumerate}[label=(W\arabic*')]
	\item For each pair $(a,b):1\le a,b\le \bbd$, \label{mu0'}$W^\varepsilon_{a,b}(dx,dk,0)$ converges (at least, a subsequence of it) to a non-negative measure $(\mu_0)_{a,b}(dx,dk)$.
\end{enumerate}

The Hamiltonian operator is
\begin{equation*}
	O^H:=\sum_{n\in\Z^{\bbd}}\alpha_n\cdot \nabla_{\beta_n}-\sum_{n,n'\in\Z^{\bbd}}\sigma(n-n')\beta_{n'}\cdot\nabla_{\alpha_n}.
\end{equation*}

\subsection{Control setting in the 1-dimensional case}
\label{Prelim1D}
We start to set up our control. First, the distance $N\in\N$ is chosen and fixed. Our controlled parameters are the $2N$ real numbers $M_N(-N),M_N(-N+1),\dots,M_N(-1),M_N(1),\dots,M_N(N)$.
We impose two assumptions on $M_N$:
\begin{enumerate}[label=($M$\arabic*)]
	\item\label{odd} $M_N(d)+M_N(-d)=0$ for $1\le |d|\le N$;
	\item\label{N/2con} $M_N(d)=0$ or $M_N(2d)=0$ for $1\le |d|\le N/2$.
\end{enumerate}

\begin{definition}
	Given a controlled $M_N$, we define the vector field
	\begin{equation}
	\lambda_{n}^{M_N}:=\sum_{d\in\Z,1\le|d|\le N}M_{N}(d)(\alpha_{n}-\alpha_{n+d})\frac{\partial}{\partial\alpha_{n-d}}+\left( \sum_{d\in\Z,1\le|d|\le N}M_{N}(d)\alpha_{n+d} \right)\frac{\partial}{\partial\alpha_n}.
	\label{vecfield}
\end{equation}

The controlled perturbation operator is defined by
\begin{equation}
	O^{M_N}:=\sum_{n\in\Z}(\lambda^{M_N}_n)^2.
	\label{OMN}
\end{equation}
\end{definition}
We will show that the operator $O^{M_N}$ conserves the total momentum and the total energy of the system in Section \ref{Conservation}.
\begin{definition}\label{1Ddefcoeff}
	Given a controlled parameter $M_N$, we define $K_{M_N}:\Z\times\Z\to \R$ as follows:
	\begin{itemize}
	\item For $d=d'=0$, 
		\begin{equation}
			K_{M_N}(0,0):= 3\sum_{d=1}^{N}M_{N}(d)^2.\label{00}
		\end{equation}
	\item For $1\le|d|\le N$ and $d'=0$,
		\begin{equation}
			K_{M_N}(d,0)= K_{M_N}(0,d):= -M_N(d)^2-\frac{1}{2}\sum_{\substack{d_1+d_2=d\\1\le|d_1|,|d_2|\le N}}M_N(d_1)M_N(d_2).\label{<=N0}
		\end{equation}
	\item For $N+1\le |d|\le 2N$ and $d'=0$,
		\begin{equation}
			K_{M_N}(d,0)= K_{M_N}(0,d):= -\frac{1}{2}\sum_{\substack{d_1+d_2=d\\1\le|d_1|,|d_2|\le N}}M_N(d_1)M_N(d_2).\label{>N0}
		\end{equation}
	\item For $0<|d|=|d'|$ and $|d|\le N/2$,
		\begin{equation}
			K_{M_N}(d,d'):= -\frac{1}{4}\sum_{\substack{d_1+d_2=d\\1\le|d_1|,|d_2|\le N}}M_N(d_1)M_N(d_2)-M_N(d)M_N(2d).\label{=<=N/2}
		\end{equation}
	\item For $0<|d|=|d'|$ and $N/2<|d|\le 2N$,
		\begin{equation}
			K_{M_N}(d,d'):= -\frac{1}{4}\sum_{\substack{d_1+d_2=d\\1\le|d_1|,|d_2|\le N}}M_N(d_1)M_N(d_2).\label{=>N/2}
		\end{equation}
	\item For $1\le |d|\le N<|d'|\le 2N$ and $|d-d'|\le N$,
		\begin{equation}
			K_{M_N}(d,d'):= \frac{1}{2}M_N(d)M_N(d'-d).\label{d<N<-}
		\end{equation}
	\item For $1\le |d|\le N<|d'|\le 2N$ and $|d+d'|\le N$,
		\begin{equation}
			K_{M_N}(d,d'):= \frac{1}{2}M_N(d)M_N(d'+d).\label{d<N<+}
		\end{equation}
	\item For $1\le |d'|\le N<|d|\le 2N$ and $|d'-d|\le N$,
		\begin{equation}
			K_{M_N}(d',d):= \frac{1}{2}M_N(d')M_N(d-d').\label{d'<N<-}
		\end{equation}
	\item For $1\le |d'|\le N<|d|\le 2N$ and $|d'+d|\le N$,
		\begin{equation}
			K_{M_N}(d',d):= \frac{1}{2}M_N(d')M_N(d+d').\label{d'<N<+}
		\end{equation}
	\item For $1\le|d|,|d'|\le N$ and $1\le |d'-d|,|d'+d|\le N$,
		\begin{equation}
			\begin{split}
				K_{M_N}(d,d')=K_{M_N}(d',d)&:= \frac{1}{2}\left( M_N(d')M_N(d-d')+M_N(d)M_N(d'-d)\right.\\
				&\hspace{4em}\left.-M_N(d')M_N(d+d')-M_N(d)M_N(d+d') \right).
			\end{split}
			\label{both<=N}
		\end{equation}
	\item For $1\le|d|,|d'|\le N$ and $1\le|d'-d|\le N<|d'+d|$,
		\begin{equation}
			K_{M_N}(d,d')=K_{M_N}(d',d):= \frac{1}{2}\left( M_N(d')M_N(d-d')+M_N(d)M_N(d'-d)\right).\label{diff<=N}
		\end{equation}
	\item For $1\le|d|,|d'|\le N$ and $1\le|d'+d|\le N<|d'-d|$,
		\begin{equation}
			K_{M_N}(d,d')=K_{M_N}(d',d):= -\frac{1}{2}\left( M_N(d')M_N(d+d')+M_N(d)M_N(d'+d)\right).\label{sum<=N}
		\end{equation}
	\item For other cases, $K_{M_N}(d,d'):=0$.
\end{itemize}
\end{definition}

\begin{definition}\label{1Ddefker}
	We define the controlled kernel corresponding to $M_N$ by
\begin{equation}\label{Ker1D}
	\hat{K}_{M_N}(k,k'):=\sum_{d,d'\in\Z}e^{-2\pi i(dk+d'k')}K_{M_N}(d,d')=\sum_{|d|,|d'|\le 2N}K_{M_N}(d,d')\cos(2\pi dk)\cos(2\pi d'k').
\end{equation}

\end{definition}

Referring to \eqref{sysmodel}, the controlled perturbation system governing the wave is
\begin{equation}
	\begin{split}
		\dot{\beta}_n(t)&= \alpha_n(t),\\
		d\alpha_{n}(t)&= -\sigma\star\beta_{n}(t)dt-2\varepsilon\left( \sum_{|d|\le 2N}K_{M_N}(d,0)\alpha_{n+d}(t) \right)dt+\sqrt{2\varepsilon}\sum_{|d|\le N}\lambda_{n+d}^{M_N}\alpha_n(t)Y_{n+d}(dt).
	\end{split}
	\label{system}
\end{equation}

	\subsection{Control setting in the multi-dimensional case}\label{Highd}
	
	\subsubsection{Notion of path}
	To define a control on multi-dimensional perturbation, we introduce the notion of path. As a path on the lattice resembles a chain, we only need to make a small modification in the computations to find the transition kernels. We first roughly introduce the controls based on a path $\gamma$.
	\begin{itemize}
		\item For a simple index expansion (see \ref{consimple}), we consider a controlled vector $M_\gamma(d)\in\R^\bbd,1\le d\le N$. Each component $[M_{\gamma}(d)]_a$ is then used to define a component of vector field $\lambda_n^{M_\gamma}$.
		\item For a dual indices expansion (see \ref{condual}), a pair of indices $(a,b)$ is chosen first, then the controls $M^{a,b}_\gamma(d),1\le d\le N$ are real numbers. After that, we also define a controlled vector field $\lambda_{n,\gamma}^{a,b}$.
	\end{itemize}

	Let us now define the notion of path. The space $\Z^{\bbd}$ is split into $2^{\bbd}$ parts using $\bbd$ hyperplanes $\{n\in\Z^{\bbd}\mid[n]_{a}=0\}$ ($a$-th component of $n\in\Z^{\bbd}$) for $a=1,\ldots,\bbd$. We map each number $p\in\{0,1,\dots,2^{\bbd}-1\}$ to one of the parts. Each number $p$ is written in base $2$, i.e. $p=p_12^{\bbd-1}+p_22^{\bbd-2}+\dots+p_{\bbd}=p_1p_2\dots p_{\bbd}|_{2}$ where $p_a\in\{0,1\}$. For each $p\in\{0,1,\dots,2^{\bbd}-1\}$, we define the set $\{n\in\Z^{\bbd}:(-1)^{p_i}[n]_a\ge0,\forall a\}$ which we will call the ``set p''. For a number $p\in\{0,1,\dots,2^{\bbd}-1\}$ there exists a unique number $p'\in\{0,1,\dots,2^{\bbd}-1\}$ whose digits in base 2 are different from the digits in base 2 of $p$ (compare in order from $p_1$ to $p_\bbd$).

	\begin{definition}
	On the set $p$ and for a number $N \in \N$, we define:
\begin{itemize}
	\item A move on the $a$-th coordinate is the map $n\to n+(-1)^{p_a}e_{a}$, where $e_a$ is the vector with $1$ at the $a$-th coordinate and $0$ elsewhere.
	\item A path is a set of points followed by making consecutive moves. We say that a path is of length $N$ if the distance between its extremities is equal to $N$. We denote by $\Gamma_{N,p}$ the family of all paths of length $N$ on the set $p$ starting from the origin. For each $\gamma\in\Gamma_{N,p}$, we denote by $\gamma_d$ the point obtained after $d$ moves ($d\ge1$).
	\item The symmetric of $\gamma$ with respect to the origin is denoted by $-\gamma$. For each $\gamma\in\Gamma_{N,p}$, we see that $-\gamma\in\Gamma_{N,p'}$.
	\end{itemize}
\end{definition}
	\begin{center}
		\begin{figure}[H]
			\centering
			\includegraphics[width=0.65\textwidth]{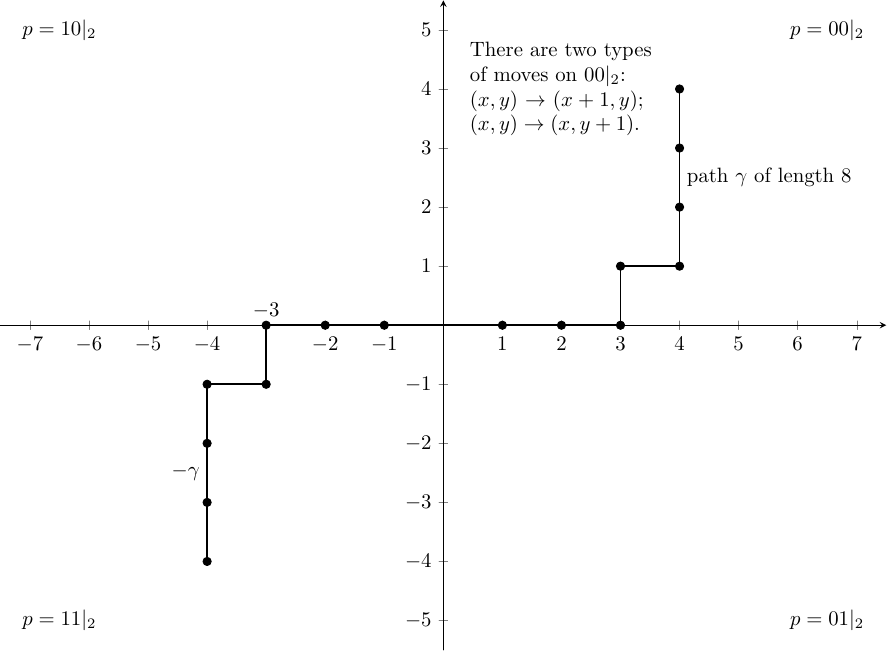}
			\caption{A path and its opposite when $\bbd=2$}
			\label{fig:path}
		\end{figure}
	\end{center}
	We note that a path that lies on the axis belongs to multiple parts. For this reason, when we mention a path we also mean that $p$ is predetermined.

	Let us establish a combinatorial lemma on the number of paths passing through a point.
	\begin{lemma}
		Given a point $D=([D]_a)_{1\le a\le\bbd}$ in the set $p\in\{0,1,\dots,2^{\bbd}-1\}$, if $\|D\|_{L^1}=|[D]_1|+\dots+|[D]_{\bbd}|\in[1,N]$ then there are \[\frac{(\|D\|_{L^1})!}{(|[D]_1|)!\dots(|[D]_{\bbd}|)!}\bbd^{N-\|D\|_{L^1}}\] paths in $\Gamma_{N,p}$ passing through $D$.
		\label{numpath}
	\end{lemma}
	\begin{proof}
		If $p_a<0$ for some $1\le a\le \bbd$ then we use the symmetry with respect to the hyperplane $\{n \in \N^{\bbd}\mid[n]_a=0\}$ and consider $\bar{p}$ with $\bar{p}_b=p_b,b\ne a$, $\bar{p}_a=-p_a$. The point $D$ has the symmetric $\bar{D}$ such that $[\bar{D}]_b=[D]_b,b\ne a$, $[\bar{D}]_a=-[D]_a$. The number of paths through $D$ in $p$ is equal to the number of paths through $\bar{D}$ in $\bar{p}$. Therefore, without loss of generality, we assume that $D$ is in the set $0$, i.e. $[D]_a\ge0$ for all $1\le a\le \bbd$.

		From the origin, we have to take $[D]_1$ moves of $+e_1$, $[D]_2$ moves of $+e_2$,\dots, $[D]_{\bbd}$ moves of $+e_\bbd$ to reach $D$. The number of paths of length $\|D\|_{L^1}$ having the end point $D$ is equal to the number of different arrangements for those moves, that is equal to
		\begin{equation*}
			\frac{(\|D\|_{L^1})!}{(|[D]_1|)!\dots(|[D]_{\bbd}|)!}.
		\end{equation*}

		For each path of length $\|D\|_{L^1}$ ending with $D$, we define a path of length $N$ by adding $N-\|D\|_{L^1}$ moves picked from $+e_1,+e_2,\dots,+e_\bbd$. There are $\bbd^{N-\|D\|_{L^1}}$ ways to pick them. The conclusion of the lemma follows.
	\end{proof}

	\subsubsection{Simple index expansion}
	\label{consimple}
	Now, we have the framework to define the control parameters $M_{\gamma}(d)\in\R^{\bbd}$. For each number $p$ and path $\gamma\in\Gamma_{N,p}$, $M_{\gamma}(d)$ is $0$ if $d< 1$ or $d>N$. Similarly to \ref{odd} and \ref{N/2con}, we impose two assumptions on the parameter $M$:
	\begin{enumerate}[label=(M\arabic*')]
		\item \label{M1'}For $p$ and $\gamma\in\Gamma_{N,p}$
	\begin{equation*}
		M_{\gamma}(d)+M_{-\gamma}(d)=0.
	\end{equation*}
\item\label{M2'} For $D\in\Z^{\bbd}$ and $1\le \|D\|_{L^1}\le N/2$,
	\begin{equation*}
		\sum_{p=0}^{2^{\bbd}-1}\sum_{\substack{\gamma \in \Gamma_{N,p}\\ D \in \gamma}}M_\gamma(\|D\|_{L^1})=0,\quad\text{or}\quad \sum_{p=0}^{2^{\bbd}-1}\sum_{\substack{\gamma \in \Gamma_{N,p}\\ 2D \in\gamma}}M_\gamma(\|2D\|_{L^1})=0.
	\end{equation*}
	\end{enumerate}

	\begin{definition}
		Given $p,N,\gamma\in\Gamma_{N,p}$, and parameter $M$, we define the vector field
	\begin{equation*}
		\lambda_{n}^{M_{\gamma}}:=\sum_{a=1}^{\bbd}\left( \sum_{d=1}^{N}[M_{\gamma}(d)]_{a}([\alpha_n]_a-[\alpha_{n+\gamma_d}]_{a})\frac{\partial}{\partial[\alpha_{n-\gamma_d}]_{a}} \right)+\sum_{a=1}^{\bbd}\left( \sum_{d=1}^{N}[M_{\gamma}(d)]_{a}[\alpha_{n+\gamma_d}]_{a}\frac{\partial}{\partial[\alpha_{n}]_{a}} \right).
	\end{equation*}
	Then, the vector field of length $N$ is defined by
	\begin{equation*}
		\begin{split}
			\lambda^{M}_{n}:= \sum_{\gamma}\lambda^{M_{\gamma}}_{n},
		\end{split}
	\end{equation*}
	where we write $\sum_{p=0}^{2^{\bbd}-1}\sum_{\gamma\in\Gamma_{N,p}}$ in short $\sum_{\gamma}$. We define the $a$-th component of $\lambda_n^{M_\gamma}$ by
	\begin{equation*}
		[\lambda^{M_\gamma}_n]_{a}:=\sum_{d=1}^{N}[M_{\gamma}(d)]_{a}([\alpha_n]_a-[\alpha_{n+\gamma_d}]_{a})\frac{\partial}{\partial[\alpha_{n-\gamma_d}]_{a}} + \sum_{d=1}^{N}[M_{\gamma}(d)]_{a}[\alpha_{n+\gamma_d}]_{a}\frac{\partial}{\partial[\alpha_{n}]_{a}}.
	\end{equation*}
The perturbed operator is
\begin{equation}
	O^{M}_s:=\sum_{n}(\lambda^{M}_n)^{2}.
	\label{OM}
\end{equation}
	\end{definition}

	\begin{definition}\label{simdefcoeff}
		For $N \in \N$, a parameter $M$ and $(a,b):1\le a,b\le \bbd$, we define the kernel coefficients $[K_M(D,D')]_{b,a}$ as follows:

\begin{itemize}
	\item $D=D'=0$,
		\begin{equation*}
			[K_{M}(0,0)]_{b,a}:=\frac{3}{2}\sum_{\substack{\gamma^1,\gamma^2,d_1,d_2\\\gamma^{1}_{d_1}=\gamma^2_{d_2}}}[M_{\gamma^1}(d_1)]_{a}[M_{\gamma^2}(d_2)]_{b}.
		\end{equation*}
	\item For $1\le\|D\|_{L^1}\le N$,
		\begin{equation*}
			\begin{split}
				[K_{M}(0,D)]_{b,a}=[K_{M}(D,0)]_{b,a}:=&-\sum_{\substack{\gamma^1,\gamma^2\\\gamma^{1}_{\|D\|_{L^1}}=\gamma^2_{\|D\|_{L^1}}=D}}[M_{\gamma^1}(\|D\|_{L^1})]_{a}[M_{\gamma^2}(\|D\|_{L^1})]_{b}\\
				&-\frac{1}{2}\sum_{\substack{\gamma^1,\gamma^2,d_1,d_2\\\gamma^{1}_{d_1}+\gamma^2_{d_2}=D}}[M_{\gamma^1}(d_1)]_{a}[M_{\gamma^2}(d_2)]_{b}.
			\end{split}
		\end{equation*}
	\item For $N+1\le\|D\|_{L^1}\le N$,
		\begin{equation*}
			[K_{M}(0,D)]_{b,a}=[K_{M}(D,0)]_{b,a}:=-\frac{1}{2}\sum_{\substack{\gamma^1,\gamma^2,d_1,d_2\\\gamma^{1}_{d_1}+\gamma^2_{d_2}=D}}[M_{\gamma^1}(d_1)]_{a}[M_{\gamma^2}(d_2)]_{b}.
		\end{equation*}
	\item For $1\le\|D\|_{L^1}\le N$,
		\begin{equation*}
			[K_{M}(D,D)]_{b,a}:=-\sum_{\substack{\gamma^1,\gamma^2\\\gamma^{1}_{\|D\|_{L^1}}=\gamma^2_{\|D\|_{L^1}}=D}}[M_{\gamma^1}(\|D\|_{L^1})]_{a}[M_{\gamma^2}(\|D\|_{L^1})]_{b}-\frac{1}{2}\sum_{\substack{\gamma^1,\gamma^2,d_1,d_2\\\gamma^{1}_{d_1}+\gamma^2_{d_2}=D}}[M_{\gamma^1}(d_1)]_{a}[M_{\gamma^2}(d_2)]_{b}.
		\end{equation*}
	\item For $1\le\|D\|_{L^1}\le N/2$,
		\begin{equation*}
			\begin{split}
				[K_{M}(D,-D)]_{b,a}:=& \sum_{\substack{\gamma^1,\gamma^2\\\gamma^{1}_{\|D\|_{L^1}}=\gamma^2_{\|D\|_{L^1}}=D}}[M_{\gamma^1}(\|D\|_{L^1})]_{a}[M_{\gamma^2}(\|D\|_{L^1})]_{b}\\
				&-\sum_{\substack{\gamma^1,\gamma^2\\\gamma^{1}_{\|D\|_{L^1}=D,\gamma^2_{\|2D\|_{L^1}}=2D}}}[M_{\gamma^1}(\|D\|_{L^1})]_{a}[M_{\gamma^2}(\|2D\|_{L^2})]_{b}\\
				&-\sum_{\substack{\gamma^1,\gamma^2\\\gamma^{1}_{\|2D\|_{L^1}=2D,\gamma^2_{\|D\|_{L^1}}=D}}}[M_{\gamma^1}(\|2D\|_{L^1})]_{a}[M_{\gamma^2}(\|D\|_{L^2})]_{b}.
			\end{split}
		\end{equation*}
	\item For $N/2<\|D\|_{L^1}\le N$,
		\begin{equation*}
			[K_{M}(D,-D)]_{b,a}:=\sum_{\substack{\gamma^1,\gamma^2\\\gamma^{1}_{\|D\|_{L^1}}=\gamma^2_{\|D\|_{L^1}}=D}}[M_{\gamma^1}(\|D\|_{L^1})]_{a}[M_{\gamma^2}(\|D\|_{L^1})]_{b}.
		\end{equation*}
	\item For $N+1\le\|D\|_{L^1}\le 2N$,
		\begin{equation*}
			[K_{M}(D,D)]_{b,a}:=-\frac{1}{2}\sum_{\substack{\gamma^1,\gamma^2,d_1,d_2\\\gamma^{1}_{d_1}+\gamma^2_{d_2}=D}}[M_{\gamma^1}(d_1)]_{a}[M_{\gamma^2}(d_2)]_{b}.
		\end{equation*}
	\item For $1\le\|D\|_{L^1}\le N<\|D'\|_{L^1}\le 2N$ and $1\le\|D-D'\|_{L^1}\le N$
		\begin{equation*}
			\begin{split}
				[K_{M}(D',D)]_{b,a}=[K_{M}(D,D')]_{b,a}:=& \frac{1}{2}\sum_{\substack{\gamma^1,\gamma^2\\\gamma^{1}_{\|D\|_{L^1}}=D\\\gamma^2_{\|D'-D\|_{L^1}}=D'-D}}[M_{\gamma^1}(\|D\|_{L^1})]_{a}[M_{\gamma^2}(\|D'-D\|_{L^1})]_{b}\\
				&+\frac{1}{2}\sum_{\substack{\gamma^1,\gamma^2\\\gamma^{1}_{\|D'-D\|_{L^1}}=D'-D\\\gamma^2_{\|D\|_{L^1}}=D}}[M_{\gamma^1}(\|D'-D\|_{L^1})]_{a}[M_{\gamma^2}(\|D\|_{L^1})]_{b}.
			\end{split}
		\end{equation*}
	\item For $1\le \|D\|_{L^1},\|D'\|_{L^1}\le N$ and $D\ne \pm D'$
		\begin{equation*}
			\begin{split}
				[K_{M}(D,D')]_{b,a}:=&\frac{1}{2}\sum_{\substack{\gamma^1,\gamma^2\\\gamma^{1}_{\|D\|_{L^1}}=D\\\gamma^2_{\|D'-D\|_{L^1}}=D'-D}}[M_{\gamma^1}(\|D\|_{L^1})]_{a}[M_{\gamma^2}(\|D'-D\|_{L^1})]_{b}\\
				&+\frac{1}{2}\sum_{\substack{\gamma^1,\gamma^2\\\gamma^{1}_{\|D'-D\|_{L^1}}=D'-D\\\gamma^2_{\|D\|_{L^1}}=D}}[M_{\gamma^1}(\|D'-D\|_{L^1})]_{a}[M_{\gamma^2}(\|D\|_{L^1})]_{b}\\
				&+\frac{1}{2}\sum_{\substack{\gamma^1,\gamma^2\\\gamma^{1}_{\|D'\|_{L^1}}=D'\\\gamma^2_{\|D-D'\|_{L^1}}=D-D'}}[M_{\gamma^1}(\|D'\|_{L^1})]_{a}[M_{\gamma^2}(\|D-D'\|_{L^1})]_{b}\\
				&+\frac{1}{2}\sum_{\substack{\gamma^1,\gamma^2\\\gamma^{1}_{\|D-D'\|_{L^1}}=D-D'\\\gamma^2_{\|D'\|_{L^1}}=D'}}[M_{\gamma^1}(\|D-D'\|_{L^1})]_{a}[M_{\gamma^2}(\|D'\|_{L^1})]_{b}\\
				&-\frac{1}{2}\sum_{\substack{\gamma^1,\gamma^2\\\gamma^{1}_{\|D'\|_{L^1}}=D'\\\gamma^2_{\|D\|_{L^1}}=D}}[M_{\gamma^1}(\|D'\|_{L^1})]_{a}[M_{\gamma^2}(\|D\|_{L^1})]_{b}\\
				&-\frac{1}{2}\sum_{\substack{\gamma^1,\gamma^2\\\gamma^{1}_{\|D\|_{L^1}}=D\\\gamma^2_{\|D'\|_{L^1}}=D'}}[M_{\gamma^1}(\|D\|_{L^1})]_{a}[M_{\gamma^2}(\|D'\|_{L^1})]_{b}.
			\end{split}
		\end{equation*}
	\item $K_{M}(D,D'):=0$ for other cases.
\end{itemize}
	\end{definition}

	\begin{definition}\label{simdefker}
		For a parameter $M$, our kernel is the matrix of coefficients
\begin{equation}\label{matrixker}
	\hat{K}_M(k,k'):=\left( \sum_{D,D'}[K_M(D,D')]_{a,b}\cos[2\pi (D\cdot k+D'\cdot k')] \right).
\end{equation}
	
	\end{definition}
This kernel corresponds to the system
\begin{equation}
	\begin{split}
		d\beta_n(t)=&\ \alpha_n(t)dt,\\
		d\alpha_n(t)=& -\sigma\star\beta_n(t)dt-2\varepsilon\left( \sum_{\|D\|_{L^1}\le N}\diag([K_{M}(D,0)]_{a,a})\alpha_{n+D} \right)dt\\
		&+\sum_{a=1}^{\bbd}\sum_{\|D\|_{L^1}\le N}\lambda^{M}_{n+D}\alpha_nY^{a}_{n+D}(dt),
	\end{split}
	\label{systemsimple}
\end{equation}
where the $Y^{a}_{n+D}$'s are independent Wiener processes.

\subsubsection{Dual indices expansion}
\label{condual}
Let us consider the energy exchange between component $a$ and component $b$ between particles on a single path. For $a,b=1,\dots,\bbd$ and $\gamma\in\Gamma_{N,p}$, we design a new control parameter $M_{\gamma}^{a,b}$.
We impose two assumptions on the controlled parameter $M^{a,b}_{\gamma}$:
\begin{enumerate}[label=(M\arabic*'')]
	\item \label{M1''}For a pair $a,b$ and a path $\gamma$
\begin{equation*}
	M^{a,b}_\gamma=M^{a,b}_{-\gamma} \text{ or }M^{a,b}_\gamma=-M^{a,b}_{-\gamma}.
\end{equation*}
\item \label{M2''}For a pair $a,b$ and a path $\gamma$
\begin{equation*}
	\sum_{1\le d_1<d_2\le N}M^{a,b}_{\gamma}(d_1)M^{a,b}_{\gamma}(d_2)=0.
\end{equation*}
\end{enumerate}
\begin{definition}
	For each $M^{a,b}$, $\gamma\in\Gamma_{N,p}$, we define the vector field
\begin{equation*}
	\begin{split}
		\lambda^{a,b}_{n,\gamma}:=&\sum_{d=1}^{N}M_{\gamma}^{a,b}(d)\left( [\alpha_{n+\gamma_{d}}]_{b}-[\alpha_{n}]_{b} \right)\sum_{d=1}^{N}M_{\gamma}^{a,b}(d)\left( \frac{\partial}{\partial[\alpha_{n+\gamma_{d}}]_{a}}-\frac{\partial}{\partial[\alpha_{n}]_{b}} \right)\\
		&-\sum_{d=1}^{N}M_{\gamma}^{a,b}(d)\left( [\alpha_{n+\gamma_{d}}]_{a}-[\alpha_{n}]_{a} \right)\sum_{d=1}^{N}M_{\gamma}^{a,b}(d)\left( \frac{\partial}{\partial[\alpha_{n+\gamma_{d}}]_{b}}-\frac{\partial}{\partial[\alpha_{n}]_{b}} \right).
	\end{split}
\end{equation*}
Using this vector field, the perturbed operator is
\begin{equation}
	O^{M}_d:=\sum_{n\in\Z^{\bbd}}\sum_{\gamma}\sum_{a,b}(\lambda^{a,b}_{n,\gamma})^2.
	\label{OMij}
\end{equation}

\end{definition}

In the case $a=b$, it is clear that $\lambda_{n,\gamma}^{a,a}=0$. Thus, we only consider $a\ne b$.
\begin{definition}
	For $a\ne b$, the coefficient $K^{a,b}$ is as follows:
	\begin{itemize}\label{dualdefcoeff}
	\item For $D=D'=0$,
		\begin{equation*}
				K^{a,b}(0,0):=\frac{1}{2}\sum_{\gamma}\left( \left( \sum_{d=1}^{N}M_{\gamma}^{a,b}(d) \right)^{2}+\left( \sum_{d=1}^{N}(M^{a,b}_{\gamma}(d))^2 \right) \right)^{2}.
		\end{equation*}
	\item For $1\le \|D\|_{L^1}\le N$,
		\begin{equation*}
			\begin{split}
				K^{a,b}(0,D)=K^{a,b}(D,0):=&-\sum_{\substack{\gamma\\\gamma_{\|D\|_{L^1}}=D}}M^{a,b}_{\gamma}(\|D\|_{L^1})\left(\sum_{d=1}^{N}M^{a,b}_{\gamma}(d)\right)\left( \left( \sum_{d=1}^{N}M_{\gamma}^{a,b}(d) \right)^{2}+\left( \sum_{d=1}^{N}(M^{a,b}_{\gamma}(d))^2 \right) \right)\\
				&+\sum_{\substack{\gamma,d_1,d_2\\\gamma_{d_2}-\gamma_{d_1}=D}}M^{a,b}_{\gamma}(d_1)M^{a,b}_{\gamma}(d_2)\left( \left( \sum_{d=1}^{N}M_{\gamma}^{a,b}(d) \right)^{2}+\left( \sum_{d=1}^{N}(M^{a,b}_{\gamma}(d))^2 \right) \right).
			\end{split}
		\end{equation*}
	\item For $1\le\|D\|_{L^1},\|D'\|_{L^1}\le N$ and $D,D'$ are both on a same set $p$ (or $D$ on set $p$ and $D'$ on set $p'$),
		\begin{equation*}
			\begin{split}
				K^{a,b}(D,D'):=&\sum_{\substack{\gamma\\\gamma_{\|D\|_{L^1}=D}\\\gamma_{\|D'\|_{L^1}}=D'}}\left( \sum_{d=1}^{N}M_{\gamma}^{a,b}(d) \right)^{2}M^{a,b}_{\gamma}(\|D\|_{L^1})M^{a,b}_{\gamma}(\|D'\|_{L^1})\\
				&-\sum_{\substack{\gamma,d_1,d_2\\\gamma_{\|D\|_{L^1}}=D\\\gamma_{d_2}-\gamma_{d_1}=D'}}\left(M^{a,b}_{\gamma}(\|D\|_{L^1})\right)\left(\sum_{d=1}^{N}M^{a,b}_{\gamma}(d)\right)M^{a,b}_{\gamma}(d_1)M^{a,b}_{\gamma}(d_2).
			\end{split}
		\end{equation*}
		We change $\gamma$ to $-\gamma$ for $M^{a,b}_{\gamma}$ in case of $D$ on set $p$ and $D'$ on set $p'$.
	\item $K^{a,b}(D,D'):=0$ for other cases.
\end{itemize}
\end{definition}

\begin{definition}\label{dualdefker}
	For each $a\ne b$ and $M^{a,b}$, we define the dual indices kernel
\begin{equation}\label{kerij}
	\hat{K}^{a,b}(k,k'):=\sum_{D,D'}K^{a,b}(D,D')\cos(2\pi D\cdot k)\cos(2\pi D'\cdot k').
\end{equation}

\end{definition}

The kernel $\widehat{K}^{a,b}$ corresponds to the system
\begin{equation}
	\begin{split}
		d\beta_{n}(t)=&\ \alpha_{n}(t)dt,\\
		d\alpha_{n}(t)=& -\sigma\star\beta_{n}(t)dt-2\varepsilon\left( \sum_{a,b}\sum_{\|D\|_{L^1}\le N}K^{a,b}(D,0)\alpha_{n+D}(t) \right)dt\\
		&+\sqrt{2\varepsilon}\sum_{a,b}\sum_{\gamma}\lambda^{a,b}_{n,\gamma}\alpha_{n}(t)Y^{a,b}_{n,\gamma}(dt)+\sum_{a,b}\sum_{1\le\|D\|_{L^1}\le N}\sum_{\gamma_{\|D\|_{L^1}}=-D}\lambda^{a,b}_{n+D,\gamma}\alpha_n(t)Y^{a,b}_{n+D,\gamma}(dt),
	\end{split}
	\label{systemij}
\end{equation}
where the $Y^{a,b}_{n,\gamma}$'s are independent Wiener processes.

\section{Main Results}

\subsection{The 1-dimensional control results}
We start with a lemma.
\begin{lemma}
	For the coupling $\sigma$ and $M_N$ satisfying Assumptions \ref{sig1} - \ref{sig4}, \ref{odd} and \ref{N/2con}, we have
\begin{equation*}
	\sum_{|d|\le 2N}K_{M_N}(d,d')=0
	%\label{sum=0}
\end{equation*}
for every $d'$ such that $|d'|\le 2N$. As a consequence,
	\begin{equation*}
		\sup_{k\in\T}\left|\frac{\hat{K}_{M_N}(k,k')}{\omega(k)}\right|<\infty.
	\end{equation*}
	\label{quotientbound}
\end{lemma}
Using the kernel defined in \eqref{Ker1D}, we also define the collision operator on the Schwartz space $\sS(\R\times\T)$ or on the space of continuous real-valued functions on $\T$ by
\begin{equation*}
	O^{M_N}_{col}S(x,k):=2\int_{\T}\hat{K}_{M_N}(k,k')\left( S(x,k')-S(x,k) \right)dk'.
\end{equation*}

\begin{theorem}
	Considering the wave system governed by \eqref{system}, we assume that the coupling $\sigma$, the initial state and $M_N$ satisfy Assumptions \ref{sig1} - \ref{sig4}, \ref{mean0} - \ref{boundenergy}, \ref{mu0}, \ref{odd} and \ref{N/2con}. Then, $E^\varepsilon(t)$ subsequentially weakly-* converges as $\varepsilon\to0$ to $\nu(t)=\nu(t,dk)$, which is the solution of
	\begin{equation}
\label{energytransport}
		\frac{d}{dt}\langle S,\nu(t)\rangle=\langle O^{M_N}_{col}S,\nu(t)\rangle
	\end{equation}
	for every bounded real-valued function $S$ on $\T$, for all $t \in [0,T]$, with $\nu(0,dk)=\nu_0(dk)$.
	\label{Enerlem}
	Moreover, under Assumption \ref{nopin} in the no-pinning case,
	\begin{equation}
\label{nopinex}
\lim_{R\to0}\limsup_{\varepsilon\to0}\frac{\varepsilon}{2}\int_{|k|<R}\E_{\varepsilon}\left[ |\hat{\phi}(k,t/\varepsilon)|^2 \right]dk=0, \quad\forall t\in[0,T].
	\end{equation}	
\end{theorem}

\begin{theorem}
	Considering the wave system governed by \eqref{system}, we assume that the coupling $\sigma$, the initial state and parameter $M_N$ satisfy Assumptions \ref{sig1} - \ref{sig4}, \ref{mean0} - \ref{nopin}, \ref{mu0}, \ref{odd} and \ref{N/2con}. Then, $W^{\varepsilon}(t)$ subsequentially weakly-* converges as $\varepsilon\to0$ to $\mu(t)=\mu(t,dx,dk)$, which is the solution of the Boltzmann equation
	\begin{equation*}
		\frac{d}{dt}\langle S,\mu(t)\rangle=\frac{1}{2\pi}\langle\omega'(k)\frac{\partial S}{\partial x},\mu(t)\rangle+\langle O^{M_N}_{col}S,\mu(t)\rangle
	\end{equation*}
	for any test function $S$ in the Schwartz space, for all $t\in[0,T]$, with $\mu(0,dx,dk)=\mu_0(dx,dk)$.
	\label{main}
\end{theorem}

\begin{remark}
	We control the Boltzmann equation in Theorem \ref{main} by modifying the kernel, through the action of $M_N$. Due to the constraints on $M_N$, we see that we can freely choose $M_{N}(d),N/2<d\le N$ and then choose $M_N(d)=0,1\le d\le N/2$. Therefore, we can choose $N-\lfloor N/2\rfloor$ values for $K_{M_N}$.
	\label{rmk}
\end{remark}

\subsection{The multi-dimensional control results}
We extend Theorem \ref{main} in two directions mentioned in Section \ref{Highd} using similar proofs. In the case of a simple index expansion, we define
\begin{equation*}
	\begin{split}
		[O^M_{col}S(x,k)]_{a,b}:=& \int_{\T^\bbd}([\hat{K}_M(k,k')]_{a,b}+[\hat{K}_M(-k,k')]_{a,b})S_{a,b}(x,k')dk'\\
		&-\int_{\T^\bbd}([\hat{K}_M(k,k')]_{a,a}+[\hat{K}_M(k,k')]_{b,b})S_{a,b}(x,k)dk'
	\end{split}
\end{equation*}
for $S=(S_{a,b})_{1\le a,b\le \bbd},S\in\sS(\R^\bbd\times\T^\bbd,\mathbb{M}_{\bbd})$ and $\hat{K}_M$ is defined in \eqref{matrixker}. The collision operator is defined by the matrix $O^M_{col}S(s,k)=([O^M_{col}S(x,k)]_{a,b})_{1\le a,b\le \bbd}$.
\begin{theorem}
	Considering the wave system governed by \eqref{systemsimple}, we assume that the coupling $\sigma$, the initial state and parameter $M$ satisfy Assumptions \ref{sig1'} - \ref{sig4'}, \ref{mean0'} - \ref{nopin'}, \ref{mu0'}, \ref{M1'} and \ref{M2'}. Then, the Wigner distribution $W^{\varepsilon}(t)$ subsequentially weakly-* converges as $\varepsilon\to0$ to $\mu(t)$, which is the solution of the Boltzmann equation
	\begin{equation*}
		\frac{d}{dt}\langle S,\mu(t)\rangle=\frac{1}{2\pi}\langle\nabla\omega(k)\cdot\nabla_xS,\mu(t)\rangle+\langle O^M_{col}S,\mu(t)\rangle
	\end{equation*}
	for any test function $S$ in the Schwartz space, for all $t\in[0,T]$, with $\mu(0,dx,dk)=\mu_0(dx,dk)$.
	\label{simplemain}
\end{theorem}

In the case of a dual indices expansion, for any $a,b$, we define
\begin{equation*}
	O^{a,b}_{col}S(x,k):=\begin{cases}
		\sum_{a'\ne a}\int_{\T^\bbd}\hat{K}^{a,a'}(k,k')(S_{a'}(x,k')-S_{a,a}(x,k))dk',&\quad\text{when }a= b,\\
		-\int_{\T^\bbd}(\hat{K}^{a,b}(k,k')+\hat{K}^{b,a}(k,k'))(S_{a,b}(x,k')+S_{a,b}(x,k))dk',&\quad\text{when }a\ne b,
	\end{cases}
\end{equation*}
where $\hat{K}^{a,b}$ is defined by \eqref{kerij}.
The collision operator for a dual indices expansion is defined by the matrix
\begin{equation*}
	O_{col}S:=(O^{a,b}_{col}S)_{1\le a,b\le \bbd}.
\end{equation*}
\begin{theorem}
	Considering the wave system governed by \eqref{systemij}, we assume that the coupling $\sigma$, the initial state and parameter $M$ satisfy Assumptions \ref{sig1'} - \ref{sig4'}, \ref{mean0'} - \ref{nopin'}, \ref{mu0'}, \ref{M1''} and \ref{M2''}. Then, the Wigner distribution $W^{\varepsilon}(t)$ subsequentially weakly-* converges as $\varepsilon\to0$ to $\mu(t)$, which is the solution of the Boltzmann equation
	\begin{equation*}
		\frac{d}{dt}\langle S,\mu(t)\rangle=\frac{1}{2\pi}\langle\nabla\omega(k)\cdot\nabla_xS,\mu(t)\rangle+\langle O_{col}S,\mu(t)\rangle
	\end{equation*}
	for any test function $S$ in the Schawartz space, for all $t\in[0,T]$, with $\mu(0,dx,dk)=\mu_{0}(dx,dk)$.
	\label{dualmain}
\end{theorem}

\subsection{Achievable Target Kernels}
In the proposition below, we prove that target kernels of the form \eqref{TargetKernel}, are achievable. 

\begin{proposition}[Achievable Target Kernels]
	In the 1-dimensional case, given any $N$, there exist $P_{a,b}(x,y),1\le a\le b\le N-\lfloor N/2\rfloor$ and a perturbation in \eqref{system} yielding the kernel
	\begin{equation}\label{TargetKernel}
		\mathcal{K}(k,k')=\sin^2(\pi k)\sin^2(\pi k')v(\cos(\pi k),\cos(\pi k')),
	\end{equation}
	where $v\in V=\left\{\sum_{1\le a\le b\le N-\lfloor N/2\rfloor}C_aC_bP_{a,b}(x,y)\mid C_a\in\R\right\}$.
	\label{propd1}
\end{proposition}

In the multi-dimensional case, obtaining such a result is much more difficult and beyond the objective of the present article. On one hand, we keep the properties $K_M(D,D')=K_M(D',D)$ and $K_M(D,D')=K_M(-D,-D')$, they are enough to derive the kernel from the controlled vectors $M_{\gamma}$. On the other hand, when we pass from dimension one to higher dimension, we lose the fact that $K_M(D,D')=K_M(-D,D')$. This significantly complicates the analysis. Furthermore, the number of free parameters in the multi-dimensional case increases exponentially with respect to the dimension.

\section{Conservation of momentum and energy - Generating operator}\label{Conservation}
\subsection{The 1-dimensional case}
We claim that $\lambda_{n}^{M_N}$ conserves the local momentum and energy. Indeed, by \ref{odd} we have
\begin{equation}
\label{1Dmoment}
	\lambda_{n}^{M_N}\left(\sum_{d'=-N}^{N}\alpha_{d'}\right)=\sum_{d\in\Z,1\le|d|\le N}M_{N}(d)(\alpha_{n}-\alpha_{n+d})+\sum_{d\in\Z,1\le|d|\le N}M_{N}(d)\alpha_{n+d}=0,
\end{equation}
and
\begin{equation}
\label{1Denergy}
	\lambda_{n}^{M_N}\left(\sum_{d'=-N}^{N}\alpha_{d'}^{2}\right)=2\sum_{d\in\Z,1\le|d|\le N}M_{N}(d)(\alpha_{n}-\alpha_{n+d})\alpha_{n-d}+2\left( \sum_{d\in\Z,1\le|d|\le N}M_{N}(d)\alpha_{n+d} \right)\alpha_{n}=0.
\end{equation}

The generating operator is identified using the following lemma.
\begin{lemma}
	The operator
\begin{equation}
	O:=O^H+\varepsilon O^{M_N}
	\label{gencon}
\end{equation}
is the generator of the stochastic differential equation \eqref{system} where $O^{M_N}$ is defined in \eqref{OMN}.
\end{lemma}
\begin{proof}
According to \eqref{phidef}, the process is
\begin{equation*}
	\phi_n(t)=\frac{1}{\sqrt{2}}\left[ \begin{matrix}
					\sigma\star\beta_n(t)\\
							\alpha_n(t)
							\end{matrix}\right].
							%\label{phivec}
					\end{equation*}
				Let us derive \eqref{gencon}. We write \eqref{system} as
				\begin{equation}
					\begin{split}
							d\left[ \begin{matrix}
											\beta_{n}(t)\\ \alpha_{n}(t)
									\end{matrix}\right]&=  \left[ \begin{matrix}
												\alpha_{n}(t)\\-\sigma\star \beta_{n}(t)-2\varepsilon\left( \sum_{|d|\le 2N}K_{M_N}(d,0)\alpha_{n+d}(t) \right)
												\end{matrix}\right]dt\\
												&+\sqrt{2\varepsilon}\left[ \begin{matrix}
																0&\dots&0\\ \lambda^{M_N}_{n-N}\alpha_{n}(t)&\dots&\lambda^{M_N}_{n+N}\alpha_{n}(t)
																\end{matrix}\right]d\left[ \begin{matrix}
																				dY_{n-N}(t)\\\vdots\\dY_{n+N}(t)
																				\end{matrix}\right].	
																				\end{split}
																				\label{dpsivec}
																		\end{equation}
																	By It\^{o}'s formula, the generator of SDE \eqref{dpsivec} is
																	\begin{equation*}
																		\alpha_{n}(t)\frac{\partial}{\partial \beta_{n}}-\sigma\star \beta_{n}(t)\frac{\partial}{\partial \alpha_{n}}-2\varepsilon \left( \sum_{|d|\le 2N}K_{M_N}(d,0)\alpha_{n+d}(t) \right)\frac{\partial}{\partial \alpha_{n}}
																		+\varepsilon\left[ \sum_{d=-N}^{N}\left( \lambda_{n+d}^{M_N}\alpha_{n}(t) \right)^{2} \right]\frac{\partial^2}{\partial \alpha_{n}^2}.
																\end{equation*}

															The term $\alpha_n\frac{\partial}{\partial \beta_n}-\sigma\star\beta_n\frac{\partial}{\partial\alpha_n}$ coincides with the Hamiltonian $O^H$. Thus, it remains to show that
															\begin{equation*}
																O^{M_N}\left(f(\alpha_n)\right)= \left[ -2\left( \sum_{|d|\le2N}K_{M_N}(d,0)\alpha_{n+d}(t) \right)\frac{\partial}{\partial \alpha_{n}}+\left[ \sum_{d=-N}^{N}\left( \lambda^{M_N}_{n+d}\alpha_{n}(t) \right)^{2} \right]\frac{\partial^2}{\partial \alpha_{n}^2} \right]f(\alpha_n),\forall n\in\Z,
														\end{equation*}
													where $f$ is a test function which is twice differentiable and compactly supported. We compute 
													\begin{equation}
														\begin{split}
																O^{M_N}(f(\alpha_n))&= \sum_{1\le|d|\le N}(\lambda_{n+d}^{M_N})^2f(\alpha_n)+(\lambda_{n}^{M_N})^2f(\alpha_n)\\
																	&= \sum_{1\le|d|\le N}\lambda_{n+d}^{M_N}M_N(d)(\alpha_{n+d}-\alpha_{n+2d})\frac{\partial f(\alpha_n)}{\partial\alpha_n}+\lambda_{n}^{M_N}\left( \sum_{1\le|d|\le N}M_N(d)\alpha_{n+d} \right)\frac{\partial f(\alpha_n)}{\partial\alpha_n}\\
																		&= \sum_{1\le|d_1|,|d_2|\le N}M_N(d_1)M_N(d_2)\alpha_{n+d_1+d_2}\frac{\partial f(\alpha_n)}{\partial\alpha_n}+\sum_{1\le|d|\le N}M_N(d)^2(\alpha_{n+d}-\alpha_{n})\frac{\partial f(\alpha_n)}{\partial\alpha_n}\\
																			&-\sum_{1\le|d|\le N}M_N(d)^2(\alpha_{n}-\alpha_{n-d})\frac{\partial f(\alpha_n)}{\partial\alpha_n}+\sum_{1\le |d|\le N}M_N(d)^2(\alpha_{n+d}-\alpha_{n+2d})^2\frac{\partial^2 f(\alpha_n)}{\partial\alpha_n^2}\\
																				&+\left( \sum_{1\le|d|\le N}M_N(d)\alpha_{n+d} \right)^2\frac{\partial^2 f(\alpha_n)}{\partial\alpha_n^2}.
																				\end{split}
																				\label{comOMN}
																		\end{equation}
																	On the other hand, the definition \eqref{vecfield} gives
																	\begin{equation*}
																		\begin{split}
																				\lambda_{n+d}^{M_N}\alpha_n&= M_{N}(d)(\alpha_{n+d}-\alpha_{n+2d})\text{ if }1\le|d|\le N,\\
																					\lambda_{n}^{M_N}\alpha_n&= \sum_{1\le|d|\le N}M_N(d)\alpha_{n+d}.
																					\end{split}
																			\end{equation*}
																		Thus, the coefficient of the second-order derivative in \eqref{comOMN} is $\sum_{d=-N}^{N}\left( \lambda_{n+d}^{M_N}\alpha_n(t) \right)^2$. The definitions \eqref{00}, \eqref{<=N0}, \eqref{>N0} also give us that the coefficient of the first-order derivative in \eqref{comOMN} is $-2\left( \sum_{|d|\le 2N}K_{M_N}(d,0)\alpha_{n+d} \right)$. This finally gives the system \eqref{system}.
																	\end{proof}
\subsection{The simple index case}
Consider the local momentum of and the local energy around $n\in\Z^{\bbd}$
	\begin{equation*}
		\sum_{n'\in\Z^{\bbd},\|n-n'\|_{L^1}\le N}\alpha_{n'}\quad\text{and}\quad\sum_{n'\in\Z^{\bbd},\|n-n'\|_{L^1}\le N}\|\alpha_{n'}\|_{L^{2}}^{2}.
	\end{equation*}
	For all $p$ and $\gamma\in\Gamma_{N,p}$ it is sufficient to prove that
		$\lambda^{M_{\gamma}}_{n}+\lambda^{M_{-\gamma}}_{n}$
conserves the local momentum and local energy.
We will design a $2N+1$-points set using path $\gamma$ and $-\gamma$, the set resembles a line in the 1-dimensional case. First, we define $n+\gamma$ as the $N+1$-points set obtained by translating $\gamma$ by $n$, i.e. $\{n,n+\gamma_1,\dots,n+\gamma_{N}\}$. We also define $n-\gamma=n+(-\gamma)$. The $2N+1$-points set is $(n+\gamma)\cup(n-\gamma)$.
If $n'\notin(n+\gamma)\cup(n-\gamma)$ then
\begin{equation*}
	(\lambda^{M_{\gamma}}_{n}+\lambda^{M_{-\gamma}}_{n})\alpha_{n'}=0\quad\text{and}\quad(\lambda^{M_{\gamma}}_{n}+\lambda^{M_{-\gamma}}_{n})\|\alpha_{n'}\|_{L^2}^{2}=0.
\end{equation*}
We only have to check the conservation for $2N+1$ points in $\{n\}\cup(n+\gamma)\cup(n-\gamma)$. We separate each index $a$ from $1$ to $\bbd$, we work on each index $a$ with similar computations as in \eqref{1Dmoment} and \eqref{1Denergy}.

\begin{lemma}
	The operator 
\begin{equation*}
	O_s:=O^H+\varepsilon O^M_s
\end{equation*}
is the generator of the SDE \eqref{systemsimple} where $O^M_s$ is defined by \eqref{OM}.
\end{lemma}
\begin{proof}
This is sufficient to check that
\begin{equation*}
	\begin{split}
		O^M_sf(\alpha_n)=&-2\sum_{a=1}^{\bbd}\sum_{\|D\|_{L^1}\le N}[K_M(D,0)]_{a,a}[\alpha_{n+D}]_a\frac{\partial f}{\partial[\alpha_n]_a}\\
				&+\sum_{a,b}\sum_{\|D\|_{L^1}\le N}([\lambda^M_{n+D}]_a[\alpha_n]_a)([\lambda^M_{n+D}]_b[\alpha_n]_b)\frac{\partial^2 f}{\partial[\alpha_n]_a\partial[\alpha_n]_b}.
				\end{split}
		\end{equation*}
	It is proved by direct computations. We see that
	\begin{align}
			O^M_sf(\alpha_n)=& \sum_{1\le\|D\|_{L^1}\le N}\left(\sum_{\gamma}\lambda^{M}_{n+D}\right)^2f(\alpha_n)+\left( \sum_{\gamma}\lambda^{M}_n \right)^{2}f(\alpha_n)\nonumber\\
				=& \sum_{1\le\|D\|_{L^1}\le N}\lambda^{M}_{n+D}\sum_{\gamma_{\|D\|_{L^1}}=D}\sum_{a=1}^{\bbd}[M_{\gamma}(\|D\|_{L^1})]_{a}([\alpha_{n+D}]_{a}-[\alpha_{n+2D}]_a)\frac{\partial f}{\partial[\alpha_n]_a}\nonumber\\
					&+\lambda^{M}_n\sum_{\gamma}\sum_{a=1}^{\bbd}\sum_{d=1}^{N}[M_\gamma(d)]_{a}[\alpha_{n+\gamma_d}]_{a}\frac{\partial f}{\partial[\alpha_n]_{a}}\nonumber\\
					=&\sum_{a=1}^{\bbd}\sum_{\gamma^1,\gamma^2,d_1,d_2}[M_{\gamma^1}(d_1)]_a[M_{\gamma^2}(d_2)]_a[\alpha_{n+\gamma^1_{d_1}+\gamma^2_{d_2}}]_{a}\frac{\partial f}{\partial[\alpha_n]_a}\label{1stsim}\\
					&+ \sum_{a=1}^{\bbd}\sum_{1\le\|D\|_{L^1}\le N}\sum_{\substack{\gamma^1,\gamma^2\\ \gamma^1_{\|D\|_{L^1}}=\gamma^{2}_{\|D\|_{L^1}}}=D}[M_{\gamma^1}(\|D\|_{L^1})]_{a}[M_{\gamma^2}(\|D\|_{L^1})]_{a}([\alpha_{n+D}]_a-[\alpha_{n}]_a)\frac{\partial f}{\partial[\alpha_n]_a}\label{2ndsim}\\
					&- \sum_{a=1}^{\bbd}\sum_{1\le\|D\|_{L^1}\le N}\sum_{\substack{\gamma^1,\gamma^2\\ \gamma^1_{\|D\|_{L^1}}=\gamma^{2}_{\|D\|_{L^1}}}=D}[M_{\gamma^1}(\|D\|_{L^1})]_{a}[M_{\gamma^2}(\|D\|_{L^1})]_{a}([\alpha_{n}]_a-[\alpha_{n-D}]_a)\frac{\partial f}{\partial[\alpha_n]_a}\label{3rdsim}\\
									&+\sum_{a,b}\sum_{1\le \|D\|_{L^1}\le N}\sum_{\substack{\gamma^1,\gamma^2\nonumber\\ \gamma^1_{\|D\|_{L^1}}=\gamma^2_{\|D\|_{L^1}}=D}}[M_{\gamma^1}(\|D\|_{L^1})]_a[M_{\gamma^2}(\|D\|_{L^1})]_b\\
									&\times([\alpha_{n+D}]_{a}-[\alpha_{n+2D}]_a)([\alpha_{n+D}]_{b}-[\alpha_{n+2D}]_b)\frac{\partial^2 f}{\partial[\alpha_n]_a\partial[\alpha_n]_b}\label{4thsim}\\
									&+\sum_{a,b}\left(  \sum_{\gamma^1}\sum_{d_1=1}^{N}[M_{\gamma^1}(d_1)]_{a}[\alpha_{n+\gamma^1_{d_1}}]_{a}\right)\left( \sum_{\gamma^2}\sum_{d_2=1}^{N}[M_{\gamma^2}(d_2)]_{b}[\alpha_{n+\gamma^2_{d_2}}]_{b} \right)\frac{\partial^2 f}{\partial[\alpha_n]_a\partial[\alpha_n]_b}.\label{5thsim}
									\end{align}

									The coefficient of $[\alpha_n]_a\frac{\partial f}{\partial [\alpha_n]_a}$ is computed from $\gamma_{d_1}^1+\gamma_{d_2}^2=0$ in \eqref{1stsim}, and from \eqref{2ndsim},\eqref{3rdsim}. They give $-2[K_{M}(0,0)]_{a,a}$. For the coefficient of $[\alpha_{n+D}]_a\frac{\partial f}{\partial [\alpha_n]_a}$, \eqref{1stsim} gives $\sum_{\substack{\gamma^1,\gamma^2,d_1,d_2\\ \gamma_{d_1}^1+\gamma_{d_2}^2=D}}[M_{\gamma^1}(d_1)]_a[M_{\gamma^2}(d_2)]_a$. Also, \eqref{2ndsim} and \eqref{3rdsim} give $2\sum_{\substack{\gamma^1,\gamma^2\\ \gamma^1_{\|D\|_{L^1}}=\gamma^{2}_{\|D\|_{L^1}}}=D}[M_{\gamma^1}(\|D\|_{L^1})]_{a}[M_{\gamma^2}(\|D\|_{L^1})]_{a}$. Thus, the coefficient is $-2[K_{M}(D,0)]_{a,a}$. The first-order derivative part is done.

									For the second-order derivative part, we fix $a,b$ and compute
									\begin{align}
										\sum_{\|D\|_{L^1}\le N}([\lambda^M_{n+D}]_a[\alpha_n]_a)([\lambda^M_{n+D}]_b[\alpha_n]_b)=& ([\lambda^M_{n}]_a[\alpha_n]_a)([\lambda^M_{n}]_b[\alpha_n]_b)\label{com2deri1}\\
										&+\sum_{1\le\|D\|_{L^1}\le N}([\lambda^M_{n+D}]_a[\alpha_n]_a)([\lambda^M_{n+D}]_b[\alpha_n]_b)\label{com2deri2}.
									\end{align}
									We can now see that \eqref{com2deri1} corresponds to the coefficient in \eqref{5thsim}, and \eqref{com2deri2} corresponds to \eqref{4thsim} by each $D$.
									\end{proof}
\subsection{The dual indices case}
We observe that the operator \eqref{OMij} conserves the total momentum and total energy. It is sufficient to check that $\lambda^{a,b}_{n,\gamma}$ conserves
\begin{equation*}
	\sum_{n'\in(n+\gamma)}[\alpha_{n'}]_{a},\quad\sum_{n'\in(n+\gamma)}[\alpha_{n'}]_{b},\quad\sum_{n'\in(n+\gamma)}([\alpha_{n'}]_{a}^{2}+[\alpha_{n'}]_{b}^{2}).
\end{equation*}
Indeed, we have
\begin{equation}
\label{dualmoment}
	\left( \sum_{d=1}^{N}M_{\gamma}^{a,b}(d)\frac{\partial}{\partial[\alpha_{n+\gamma_{d}}]_{a}}-\sum_{d=1}^{N}M_{\gamma}^{a,b}(d)\frac{\partial}{[\alpha_{n}]_{a}} \right)\sum_{n'\in\{n\}\cup(n+\gamma)}[\alpha_{n'}]_{a}=\left( \sum_{d=1}^{N}M_{\gamma}^{a,b}(d)-\sum_{d=1}^{N}M_{\gamma}^{a,b}(d) \right)=0.
\end{equation}
We prove the conservation of the $b$-th component of the moment similarly to \eqref{dualmoment}. For the energy we have
\begin{equation*}
	\begin{split}
		\lambda^{a,b}_{n,\gamma}\sum_{n'\in\{n\}\cup(n+\gamma)}&([\alpha_{n'}]_{a}^{2}+[\alpha_{n'}]_{b}^{2})\\
		=& 2\left( \sum_{d=1}^{N}M_{\gamma}^{a,b}(d)[\alpha_{n+\gamma_{d}}]_{b}-\sum_{d=1}^{N}M_{\gamma}^{a,b}(d)[\alpha_{n}]_{b} \right)\left( \sum_{d=1}^{N}M_{\gamma}^{a,b}(d)[\alpha_{n+\gamma_{d}}]_{a}-\sum_{d=1}^{N}M_{\gamma}^{a,b}(d)[\alpha_{n}]_{a} \right)\\
		&-2\left( \sum_{d=1}^{N}M_{\gamma}^{a,b}(d)[\alpha_{n+\gamma_{d}}]_{a}-\sum_{d=1}^{N}M_{\gamma}^{a,b}(d)[\alpha_{n}]_{a} \right)\left( \sum_{d=1}^{N}M_{\gamma}^{a,b}(d)[\alpha_{n+\gamma_{d}}]_{b}-\sum_{d=1}^{N}M_{\gamma}^{a,b}(d)[\alpha_{n}]_{b} \right)\\
		=& 0.
	\end{split}
\end{equation*}

\begin{lemma}
	The operator
\begin{equation*}
	O_d:=O^H+\varepsilon O^M_d
\end{equation*}
is the generator of the SDE \eqref{systemij} where $O^M_d$ is defined by \eqref{OMij}.
\end{lemma}
\begin{proof}
We need to prove that
\begin{equation*}
	\begin{split}
			O^M_d=& -2\left( \sum_{a,b}\sum_{\|D\|_{L^1}\le N}K^{a,b}(D,0)\alpha_{n+D}(t) \right)\cdot\frac{\partial}{\partial\alpha_n}\\
				&+\sum_{a,b}\sum_{\gamma}\sum_{a_1,a_2\in\{a,b\}}[(\lambda^{a,b}_{n,\gamma}\alpha_n(t))(\lambda^{a,b}_{n,\gamma}\alpha_n(t))^{\top}]_{a_1,a_2}\frac{\partial^2}{\partial[\alpha_n]_{a_1}\partial[\alpha_n]_{a_2}}\\
					&+\sum_{a,b}\sum_{1\le\|D\|_{L^1}\le N}\sum_{\gamma_{\|D\|_{L^1}=-D}}\sum_{a_1,a_2\in\{a,b\}}[(\lambda^{a,b}_{n+D,\gamma}\alpha_n(t))(\lambda^{a,b}_{n+D,\gamma}\alpha_n(t))^{\top}]_{a_1,a_2}\frac{\partial^2}{\partial[\alpha_n]_{a_1}\partial[\alpha_n]_{a_2}}
					\end{split}
			\end{equation*}
		Computing $O^M_df$, we get
		\begin{equation*}
			O^M_df(\alpha_n)= \sum_{a,b}\sum_{\gamma}\sum_{1\le\|D'\|_{L^1}\le N}(\lambda^{a,b}_{n+D',\gamma})^2f(\alpha_n)+\sum_{a,b}\sum_{\gamma}(\lambda^{a,b}_{n,\gamma})^{2}f(\alpha_n).
	\end{equation*}
	For $D'\in\Z^{\bbd},1\le \|D'\|_{L^1}\le N$, we compute
\begin{align*}
	\sum_{\gamma}(\lambda^{a,b}_{n+D',\gamma})^{2}f(\alpha_n)=& -\sum_{\gamma_{\|D'\|_{L^1}}=-D'}\sum_{d=1}^{N}M^{a,b}_{\gamma}(d)([\alpha_{n+\gamma_{d}+D'}]_{a}-[\alpha_{n+D'}]_a)\numberthis\label{nD'}\\
	&\times\left( \left( \sum_{d=1}^{N}M_{\gamma}^{a,b}(d) \right)^{2}+\left( \sum_{d=1}^{N}(M^{a,b}_{\gamma}(d))^2 \right) \right)M^{a,b}_{\gamma}(\|D'\|_{L^1})\frac{\partial f(\alpha_n)}{\partial[\alpha_n]_{a}}\\
	&- \sum_{\gamma_{\|D'\|_{L^1}}=-D'}\sum_{d=1}^{N}M^{a,b}_{\gamma}(d)([\alpha_{n+\gamma_{d}+D'}]_{b}-[\alpha_{n+D'}]_b)\\
	&\times\left( \left( \sum_{d=1}^{N}M_{\gamma}^{a,b}(d) \right)^{2}+\left( \sum_{d=1}^{N}(M^{a,b}_{\gamma}(d))^2 \right) \right)M^{a,b}_{\gamma}(\|D'\|_{L^1})\frac{\partial f(\alpha_n)}{\partial[\alpha_n]_{b}}\\
	&+ \sum_{\gamma_{\|D'\|_{L^1}}=-D'}\left( \sum_{d=1}^{N}M^{a,b}_{\gamma}([\alpha_{n+\gamma_{d}+D'}]_{b}-[\alpha_{n+D'}]_b) \right)^{2}\left( M^{a,b}_{\gamma}(\|D'\|_{L^1}) \right)^{2}\frac{\partial^2f(\alpha_n)}{\partial[\alpha_n]_{a}^2}\\
	&+ \sum_{\gamma_{\|D'\|_{L^1}}=-D'}\left( \sum_{d=1}^{N}M^{a,b}_{\gamma}([\alpha_{n+\gamma_{d}+D'}]_{a}-[\alpha_{n+D'}]_a) \right)^{2}\left( M^{a,b}_{\gamma}(\|D'\|_{L^1}) \right)^{2}\frac{\partial^2f(\alpha_n)}{\partial[\alpha_n]_{b}^2}\\
	&- 2\sum_{\gamma_{\|D'\|_{L^1}}=-D'}\left( \sum_{d=1}^{N}M^{a,b}_{\gamma}([\alpha_{n+\gamma_{d}+D'}]_{a}-[\alpha_{n+D'}]_a) \right)\left( M^{a,b}_{\gamma}(\|D'\|_{L^1}) \right)^{2}\\
	&\times\left( \sum_{d=1}^{N}M^{a,b}_{\gamma}([\alpha_{n+\gamma_{d}+D'}]_{b}-[\alpha_{n+D'}]_b) \right)\frac{\partial^2f(\alpha_n)}{\partial[\alpha_n]_{a}\partial[\alpha_n]_b}.
\end{align*}
We also have
\begin{align*}
	\sum_{\gamma}(\lambda^{a,b}_{n,\gamma})^{2}f(\alpha_n)=& \sum_{\gamma}\sum_{d=1}^{N}M^{a,b}_{\gamma}(d)([\alpha_{n+\gamma_{d}}]_{a}-[\alpha_n]_a)\numberthis\label{n0}\\
	&\times\left( \left( \sum_{d=1}^{N}M_{\gamma}^{a,b}(d) \right)^{2}+\left( \sum_{d=1}^{N}(M^{a,b}_{\gamma}(d))^2 \right) \right)\left( \sum_{d=1}^{N}M^{a,b}_{\gamma} \right)\frac{\partial f(\alpha_n)}{\partial[\alpha_n]_{a}}\\
	&+ \sum_{\gamma}\sum_{d=1}^{N}M^{a,b}_{\gamma}(d)([\alpha_{n+\gamma_{d}}]_{b}-[\alpha_n]_b)\\
	&\times\left( \left( \sum_{d=1}^{N}M_{\gamma}^{a,b}(d) \right)^{2}+\left( \sum_{d=1}^{N}(M^{a,b}_{\gamma}(d))^2 \right) \right)\left( \sum_{d=1}^{N}M^{a,b}_{\gamma} \right)\frac{\partial f(\alpha_n)}{\partial[\alpha_n]_{b}}\\
	&+ \sum_{\gamma}\left( \sum_{d=1}^{N}M^{a,b}_{\gamma}([\alpha_{n+\gamma_{d}}]_{b}-[\alpha_n]_b) \right)^{2}\left( \sum_{d=1}^{N}M^{a,b}_{\gamma}(d) \right)^{2}\frac{\partial^2f(\alpha_n)}{\partial[\alpha_n]_{a}^2}\\
	&+ \sum_{\gamma}\left( \sum_{d=1}^{N}M^{a,b}_{\gamma}([\alpha_{n+\gamma_{d}}]_{a}-[\alpha_n]_a) \right)^{2}\left( \sum_{d=1}^{N}M^{a,b}_{\gamma}(d) \right)^{2}\frac{\partial^2f(\alpha_n)}{\partial[\alpha_n]_{b}^2}\\
	&- 2\sum_{\gamma}\left( \sum_{d=1}^{N}M^{a,b}_{\gamma}([\alpha_{n+\gamma_{d}}]_{a}-[\alpha_n]_a) \right)\left( \sum_{d=1}^{N}M^{a,b}_{\gamma}(d) \right)^{2}\\
	&\times\left( \sum_{d=1}^{N}M^{a,b}_{\gamma}([\alpha_{n+\gamma_{d}}]_{b}-[\alpha_n]_b) \right)\frac{\partial^2f(\alpha_n)}{\partial[\alpha_n]_{a}\partial[\alpha_n]_b}.
\end{align*}

When $D=0$, \eqref{n0} gives
\begin{equation*}
	\sum_{\gamma}\left( \left( \sum_{d=1}^{N}M_{\gamma}^{a,b}(d) \right)^{2}+\left( \sum_{d=1}^{N}(M^{a,b}_{\gamma}(d))^2 \right) \right)\left( \sum_{d=1}^{N}M^{a,b}_{\gamma}(d) \right)^2\frac{\partial}{\partial[\alpha_n]_a}
\end{equation*}
and $\eqref{nD'}$ gives
\begin{equation*}
	\sum_{D'}\sum_{\gamma_{\|D'\|_{L^1}}=-D'}\left( \left( \sum_{d=1}^{N}M_{\gamma}^{a,b}(d) \right)^{2}+\left( \sum_{d=1}^{N}(M^{a,b}_{\gamma}(d))^2 \right) \right)\left( M^{a,b}_{\gamma}(\|D'\|_{L^1}) \right)^2\frac{\partial}{\partial[\alpha_n]_a}.
\end{equation*}
The sum of two term in case $D=0$ is $-2K^{a,b}(0,0)$. Now, we consider $1\le \|D\|_{L^1}\le N$, \eqref{n0} gives
\begin{equation*}
	\sum_{\gamma_{\|D\|_{L^1}}=D}\left( \left( \sum_{d=1}^{N}M_{\gamma}^{a,b}(d) \right)^{2}+\left( \sum_{d=1}^{N}(M^{a,b}_{\gamma}(d))^2 \right) \right)\left( \sum_{d=1}^{N}M^{a,b}_{\gamma}(d) \right)M^{a,b}_{\gamma}(\|D\|_{L^1})\frac{\partial}{\partial[\alpha_n]_a}.
\end{equation*}
Meanwhile, \eqref{nD'}, for the case $D=D'$, also gives the same term. For the case $D\ne D'$, \eqref{nD'} only gives non-zero term if there is $d_2$ such that $\gamma_{d_2}-\gamma_{\|D'\|_{L^1}}=D$. Thus, with $1\le\|D\|_{L^1}\le N$, we get the coefficient
\begin{equation*}
	2\sum_{\substack{\gamma,d_1,d_2\\\gamma_{d_2}-\gamma_{d_1}=D}}M^{i,j}_{\gamma}(d_1)M^{i,j}_{\gamma}(d_2)\left( \left( \sum_{d=1}^{N}M_{\gamma}^{i,j}(d) \right)^{2}+\left( \sum_{d=1}^{N}(M^{i,j}_{\gamma}(d))^2 \right) \right).
\end{equation*}
Hence, we get $-2K^{i,j}(D,0)$. For the second-order derivative parts, we see that
\begin{equation*}
	\begin{split}
			\lambda_{n,\gamma}^{a,b}\alpha_n=& -\sum_{d=1}^{N}M_{\gamma}^{a,b}(d)\left( [\alpha_{n+\gamma_{d}}]_{b}-[\alpha_{n}]_{b} \right)\sum_{d=1}^{N}M_{\gamma}^{a,b}(d)e_{a}\\
				&+\sum_{d=1}^{N}M_{\gamma}^{a,b}(d)\left( [\alpha_{n+\gamma_{d}}]_{a}-[\alpha_{n}]_{a} \right)\sum_{d=1}^{N}M_{\gamma}^{a,b}(d)e_{b}.
				\end{split}
		\end{equation*}
	Hence, $(\lambda^{a,b}_{n,\gamma}\alpha_n)(\lambda^{a,b}_{n,\gamma}\alpha_n)^{\top}$ matches the second-order derivative part of \eqref{n0}. Similarly, by a direct computation, $(\lambda^{a,b}_{n+D,\gamma}\alpha_n)(\lambda^{a,b}_{n+D,\gamma}\alpha_n)^{\top}$ matches the second-order derivative part of \eqref{nD'}.
\end{proof}

\section{Computation of the kernel}

\subsection{The 1-dimensional case}\label{comker}
In this section, we show how to compute the kernel $K_{M_N}$ defined by \eqref{00} and \eqref{sum<=N}. We start with the computation of the generator over the energy density $O|\hat{\phi}|^2$, where $O$ is defined in \eqref{gencon}. By \eqref{O^Hiden}, we get $O^H|\hat{\phi}|^2=0$. Thus, we focus on $O^{M_N}|\hat{\phi}|^2$. Since $O^{M_N}$ is a second-order operator we get
\begin{equation*}
	O^{M_N}|\hat{\phi}|^2=(O^{M_N}\hat{\phi})\hat{\phi}^{*}+(O^{M_N}\hat{\phi}^{*})\hat{\phi}+2\sum_{n}\left( \lambda^{M_N}_{n}\hat{\phi}(k) \right)\left( \lambda^{M_N}_n\hat{\phi}^{*}(k) \right).
\end{equation*}
The kernel will come from $\sum_{n}\left( \lambda^{M_N}_{n}\hat{\phi}(k) \right)\left( \lambda^{M_N}_n\hat{\phi}^{*}(k) \right)$.
We compute
\begin{align}
		\sum_{n\in\Z}&(\lambda_{n}^{M_N}\phi_{n'})(\lambda_{n}^{M_N}\phi_{n''}^{*})\nonumber\\
		&= \frac{1}{2}\sum_{n\in\Z}\left[\sum_{\substack{d\in\Z\\1\le|d|\le N}}M_N(d)(\alpha_{n}-\alpha_{n+d})\delta_{n-n',d}+\left( \sum_{\substack{d\in\Z\\1\le|d|\le N}}M_N(d)\alpha_{n+d} \right)\delta_{n,n'}\right]\nonumber\\
		&\times\left[ \sum_{\substack{d\in\Z\\1\le|d|\le N}}M_N(d)(\alpha_{n}-\alpha_{n+d})\delta_{n-n'',d}+\left( \sum_{\substack{d\in\Z\\1\le|d|\le N}}M_N(d)\alpha_{n+d} \right)\delta_{n,n''} \right]\nonumber\\
		&= -\frac{1}{2}\sum_{d=-2N}^{2N}\left( \sum_{\substack{1\le|d_1|,|d_2|\le N\\d_1+d_2=d'}}M_N(d_1)M_N(d_2)(\alpha_{n'-d_1}-\alpha_{n'-2d_1})(\alpha_{n'-d_1}-\alpha_{n'-d_1+d_2})\right)\delta_{n'-n'',d'}\label{sm1}\\
		& -\frac{1}{2}\sum_{1\le|d_1|\le N}M_N(d_1)(\alpha_{n'-d_1}-\alpha_{n'-2d_1})\left( \sum_{\substack{d_2\in\Z\\1\le|d_2|\le N}}M_N(d_2)\alpha_{n'-d_1+d_2} \right)\delta_{n'-n'',d_1}\label{sm2}\\
		& +\frac{1}{2}\sum_{1\le|d_1|\le N}M_N(d_1)(\alpha_{n'}-\alpha_{n'+d_1})\left( \sum_{\substack{d_2\in\Z\\1\le|d_2|\le N}}M_N(d_2)\alpha_{n'+d_2}\right)\delta_{n'-n'',d_1}\label{sm3}\\
		&+\frac{1}{2}\left( \sum_{\substack{d\in\Z\\1\le|d|\le N}}M_N(d)\alpha_{n'+d} \right)^{2}\delta_{n'-n'',0}.\label{sm4}
\end{align}
Let us write $\sum_{n,n',n''}(\lambda_{n}^{M_N}\phi_{n'})(\lambda_{n}^{M_N}\phi_{n''}^{*})$ in the form of $\sum_{n',n'',d,d'}K_{M_N}(d,d')\alpha_{n'}\alpha_{n'+d}\delta_{n'-n'',d'}$.
Taking into consideration the invariance of $n'-n''$ under the translation of $n',n''$ by a same number, we cancel $-1/2M_N(d_1)\alpha_{n'-d_1}(\sum_{d_2}M_N(d_2)\alpha_{n'-d_1+d_2})$ in \eqref{sm2} with $1/2M_N(d_1)\alpha_{n'}(\sum_{d_2}M_N(d_2)\alpha_{n'+d_2})$ in \eqref{sm3}.
We obtain the kernel $K_{M_N}$ as follows:
\begin{itemize}
	\item If $d=d'=0$, the term \eqref{sm1} gives $2\sum_{d=1}^{N}M_{N}(d)^2$($d_1=-d_2$), and the term \eqref{sm4} gives $\sum_{d=1}^{N}M_{N}(d)^2$. We get \eqref{00}.
	\item If $d'=0$ and $1\le |d|\le N$, the term \eqref{sm1} gives $-M_N(d)^2$($d_1=-d_2=d$), and the term \eqref{sm4} gives $-\frac{1}{2}\sum_{d_1+d_2=d,1\le|d_1|,|d_2|\le N}M_N(d_1)M_N(d_2)$. If $1\le|d'|\le N$ and $d=0$, the term \eqref{sm1} gives $-\frac{1}{2}\sum_{d_1+d_2=d,1\le|d_1|,|d_2|\le N}M_N(d_1)M_N(d_2)$, the term \eqref{sm2} gives $-\frac{1}{2}M_N(d)^2$($d_1=-d_2=d$), and the term \eqref{sm3} gives $-\frac{1}{2}M_N(d)^2$($d_1=d_2=d$). We get \eqref{<=N0}.
	\item If $d'=0$ and $N+1\le |d|\le 2N$, there is no contribution from the terms \eqref{sm1}, \eqref{sm2}, \eqref{sm3}, and the term \eqref{sm4} gives \[K_N(d,0)=-\frac{1}{2}\sum_{\substack{d_1+d_2=d\\1\le|d_1|,|d_2|\le N}}M_N(d_1)M_N(d_2).\]
		If $N+1\le|d'|\le 2N$ and $d=0$, the term \eqref{sm1} gives	\[K_N(0,d')=-\frac{1}{2}\sum_{\substack{d_1+d_2=d'\\1\le|d_1|,|d_2|\le N}}M_N(d_1)M_N(d_2).\]
Thus, we get \eqref{>N0}.
	\item If $|d'|=|d|$, the term \eqref{sm1} gives $-\frac{1}{2}\sum_{d_1+d_2=d',1\le|d_1|,|d_2|\le N}M_N(d_1)M_N(d_2)$, it also gives $-M_N(d')M_N(2d')$ ($d_1=-d',d_2=2d'$ or $d_1=2d',d_2=-d'$ in the case $|d|\le N/2$). When $|d|<N/2$, the term \eqref{sm2} gives $-1/2M_N(d')M_N(2d')$($d_1=d',d_2=-2d'$), and the term \eqref{sm3} gives $-1/2M_N(d')M_N(2d')$($d_1=d',d_2=2d'$). Therefore,
		\begin{equation*}
			\begin{split}
				K_{M_N}(d,d')&= -\frac{1}{4}\sum_{\substack{d_1+d_2=d\\1\le|d_1|,|d_2|\le N}}M_N(d_1)M_N(d_2)-M_N(d)M_N(2d)\quad\text{if }|d|\le N/2,\\
				K_{M_N}(d,d')&= -\frac{1}{4}\sum_{\substack{d_1+d_2=d\\1\le|d_1|,|d_2|\le N}}M_N(d_1)M_N(d_2)\quad\text{for other cases.}
			\end{split}
		\end{equation*}
		There is only one term $\alpha_{n'-2d_1}\alpha_{n'-d_1+d_2}$ that has the form $\alpha_n\alpha_{n+d'}$, we use the factor $1/4$ to ensure the symmetry on $K_{M_N}$. We get \eqref{=<=N/2} and \eqref{=>N/2}.
	\item If $1\le|d|\le N<|d'|\le 2N$ and $|d'|-|d|\le N$, then
		\begin{equation*}
			\begin{split}
				K_{M_N}(d,d')&= \frac{1}{2}M_N(d)M_N(d'-d)\qquad\text{if }|d'-d|\le N,\\
				K_{M_N}(d,d')&= -\frac{1}{2}M_N(d)M_N(d'+d)\qquad\text{if }|d'+d|\le N.\\
			\end{split}
		\end{equation*}
		It corresponds to the case $d_1=d,d_2=d'-d$ and $d_1=d'-d,d_2=d$ in \eqref{sm1}. Here, we also use the factor $1/2$ to ensure the symmetry of $K_{M_N}$. We get \eqref{d<N<-} and \eqref{d<N<+}.
	\item If $1\le|d'|\le N<|d|\le 2N$ and $|d|-|d'|\le N$, then
		\begin{equation*}
			\begin{split}
				K_{M_N}(d,d')&= \frac{1}{2}M_N(d')M_N(d-d')\qquad\text{if }|d-d'|\le N,\\
				K_{M_N}(d,d')&= -\frac{1}{2}M_N(d')M_N(d+d')\qquad\text{if }|d+d'|\le N.\\
			\end{split}
		\end{equation*}
		It corresponds to the case $d_1=d',d_2=d-d'$ in \eqref{sm2} and $d_1=d',d_2=d'+d$ in \eqref{sm3}. We get \eqref{d'<N<-} and \eqref{d'<N<+}.
	\item If $1\le|d'|,|d|\le N$ and $1\le|d-d'|,|d+d'|\le N$, we compute from \eqref{sm1} ($d_1=d,d_2=d'-d$ or $d_1=d-d',d_2=d$, and ensure the symmetry by the factor $1/2$), \eqref{sm2} ($d_1=d',d_2=d-d'$), \eqref{sm3} ($d_1=d',d_2=d+d'$), the coefficient is
		\[K_N(d,d')=\frac{1}{2}(M_N(d')M_N(d-d')-M_N(d')M_N(d+d')-M_N(d)M_N(d-d')-M_N(d)M_N(d+d')).\]
	\item If $1\le|d'|,|d|\le N$ and $1\le|d-d'|\le N<|d'+d|$, we compute from \eqref{sm1} ($d_1=d,d_2=d'-d$ or $d_1=d-d',d_2=d$, and ensure the symmetry by the factor $1/2$), \eqref{sm2} ($d_1=d',d_2=d-d'$), the coefficient is
		\[K_N(d,d')=\frac{1}{2}(M_N(d')M_N(d-d')-M_N(d)M_N(d-d')).\]
	\item If $1\le|d'|,|d|\le N$ and $1\le|d+d'|\le N<|d-d'|$, we compute from \eqref{sm1} ($-1/2M_N(d')M_N(d+d')$ from the previous cases), and \eqref{sm3} ($d_1=d',d_2=d+d'$), the coefficient is
		\[K_N(d,d')=\frac{1}{2}(-M_N(d')M_N(d+d')-M_N(d)M_N(d+d')).\]
	\item $K_{N}(d,d')=0$ in all other cases.
\end{itemize}
We get the coefficients $K_{M_N}$ in Definition \ref{1Ddefcoeff}.

On the kernel $\widehat{K}_{M_N}$ in Definition \ref{1Ddefker}, we have
\begin{equation*}
	\begin{split}
		\sum_{n\in\Z}\left( \lambda^{M_N}_n\hat{\phi}(k) \right)\left( \lambda^{M_N}_n\hat{\phi}^{*}(k) \right)&= \sum_{n',n''\in\Z}e^{2\pi ik(n''-n')}\sum_{n}\left( \lambda^{M_N}_n\phi_{n'} \right)\left( \lambda^{M_N}_n\phi_{n''} \right)\\
		&= \sum_{n'd,d'}K_{M_N}(-d',d)e^{-2\pi ikd}\alpha_{n'}\alpha_{n'-d'}\\
		&= \int_{\T}\sum_{d,d'}e^{-2\pi(kd+k'd')}K_{M_N}(d,d')|\hat{\alpha}(k')|^2dk'\\
		&= \int_{\T}\hat{K}_{M_N}(k,k')|\hat{\alpha}(k')|^2dk'\\
		&= \int_{\T}\hat{K}_{M_N}(k,k')\left( |\hat{\phi}(k')|^2-\frac{1}{2}\left( \hat{\phi}(k')\hat{\phi}(-k')+\hat{\phi}^{*}(k')\hat{\phi}^{*}(-k') \right) \right).
	\end{split}
\end{equation*}
We also compute
\begin{equation*}
	\begin{split}
		O^{M_N}\hat{\phi}(k)&= \frac{1}{2}\sum_{n\in\Z}e^{-2\pi ikn}\left( -2\sum_{d\in\Z}K_{M_N}(d,0)\left( \phi_{n-d}-\phi_{n-d}^{*} \right) \right)\\
		&= -\sum_{d\in\Z}e^{-2\pi ikd}K_{M_N}(d,0)\left( \hat{\phi}(k)-\hat{\phi}^{*}(-k) \right)\\
		&= -\int_{\T}\hat{K}_{M_N}(k,k')\left( \hat{\phi}(k)-\hat{\phi}^{*}(-k) \right)dk'.
	\end{split}
\end{equation*}
Thus, we get
\begin{equation*}
	(O^{M_N}\hat{\phi}(k))\hat{\phi}^{*}(k)+(O^{M_N}\hat{\phi}^{*}(k))\hat{\phi}(k)=\int_{\T}\hat{K}_{M_N}(k,k')\left( -2|\hat{\phi}(k)|^2+\left( \hat{\phi}(k)\hat{\phi}(-k)+\hat{\phi}^{*}(k)\hat{\phi}^{*}(-k) \right) \right)dk'.
\end{equation*}
We conclude that
\begin{equation*}
	O|\hat{\phi}(k)|^2=\varepsilon O^{M_N}_{col}\left( |\hat{\phi}(k)|^2-\frac{1}{2}\left( \hat{\phi}(k)\hat{\phi}(-k)+\hat{\phi}^{*}(k)\hat{\phi}^{*}(-k) \right) \right).
\end{equation*}

\subsection{The simple index case} We follow the steps in the 1-dimensional case to compute the kernel. From the parameter $M$, we derive the controlled kernel coefficients by computing
\begin{equation*}
	\sum_{n,n',n''\in\Z^{\bbd}}\sum_{\gamma^1}\sum_{\gamma^2}[\lambda^{M_{\gamma^1}}_{n}]_{a}[\phi_{n'}]_{a}[\lambda^{M_{\gamma^2}}_n]_{b}[\phi_{n''}^{*}]_{b}.
\end{equation*}

Omitting temporarily the sum over $n',n''$ in the formulation (and still using the invariance of $n'-n''$ under translation), we compute
\begin{align}
	&\sum_{n\in\Z^{\bbd}}\sum_{\gamma^1}\sum_{\gamma^2}[\lambda^{M_{\gamma^1}}_{n}]_{a}[\phi_{n'}]_{a}[\lambda^{M_{\gamma^2}}_n]_{b}[\phi^{*}_{n''}]_{b}\nonumber\\
		&=  \frac{1}{2}\sum_{n\in\Z^{\bbd}}\sum_{\gamma^1}\sum_{\gamma^2}\left[ \sum_{d=1}^{N}[M_{\gamma^1}(d)]_{a}([\alpha_{n}]_{a}-[\alpha_{n+\gamma^1_d}]_{a})\delta_{n-n',\gamma^1_d}+\left( \sum_{d=1}^{N}[M_{\gamma^1}(d)]_{a}[\alpha_{n+\gamma^1_d}]_a \right)\delta_{n,n'} \right]\nonumber\\
		&\hspace{2em}\times\left[ \sum_{d=1}^{N}[M_{\gamma^2}(d)]_{b}([\alpha_{n}]_{b}-[\alpha_{n+\gamma^2_d}]_{b})\delta_{n-n'',\gamma^2_d}+\left( \sum_{d=1}^{N}[M_{\gamma^2}(d)]_{b}[\alpha_{n+\gamma^2_d}]_b \right)\delta_{n,n''} \right]\nonumber\\
		&= -\frac{1}{2}\sum_{\substack{D'\in\Z^{\bbd}\\\|D'\|_{L^1}\le 2N}}\left( \sum_{\substack{\gamma^1,\gamma^2,d_1,d_2\\\gamma^{1}_{d_1}+\gamma^2_{d_2}=D'}}[M_{\gamma^1}(d_1)]_{a}[M_{\gamma^2}(d_2)]_{b}([\alpha_{n'-\gamma^1_{d_1}}]_{a}-[\alpha_{n'-2\gamma^1_{d_1}}]_{a})\right.\nonumber\\
		&\hspace{10em}\times\left.\vphantom{\sum_{\substack{\gamma^1,\gamma^2,d_1,d_2\\\gamma^{1}_{d_1}+\gamma^2_{d_2}=D}}}([\alpha_{n'-\gamma^1_{d_1}}]_b-[\alpha_{n'-\gamma^1_{d_1}+\gamma^2_{d_2}}]_b) \right)\delta_{n'-n'',D'}\label{sim1}\\
		&-\frac{1}{2}\sum_{\gamma^1,d_1}[M_{\gamma^1}(d_1)]_{a}([\alpha_{n'-\gamma^1_{d_1}}]_{a}-[\alpha_{n'-2\gamma^1_{d_1}}]_{a})\left( \sum_{\gamma^2,d_2}[M_{\gamma^2}(d_2)]_{b}[\alpha_{n'-\gamma^1_{d_1}+\gamma^2_{d_2}}]_b \right)\delta_{n'-n'',\gamma^1_{d_1}}\label{sim2}\\
		&+\frac{1}{2}\sum_{\gamma^2,d_2}[M_{\gamma^2}(d_2)]_{b}([\alpha_{n'}]_{b}-[\alpha_{n'+\gamma^2_{d_2}}]_{b})\left( \sum_{\gamma^1,d_1}[M_{\gamma^1}(d_1)]_{a}[\alpha_{n'+\gamma^1_{d_1}}]_a \right)\delta_{n'-n'',\gamma^2_{d_2}}\label{sim3}\\
		&+\frac{1}{2}\left( \sum_{\gamma^1}\sum_{d=1}^{N}[M_{\gamma^1}(d)]_{a}[\alpha_{n+\gamma^1_d}]_a \right)\left( \sum_{\gamma^2}\sum_{d=1}^{N}[M_{\gamma^2}(d)]_{b}[\alpha_{n+\gamma^2_d}]_b \right)\delta_{n'-n'',0}.\label{sim4}
\end{align}
We expect to write the sum in the form $\sum_{D,D'}[K_{M}(D,D')]_{b,a}[\alpha_{n'}]_{a}[\alpha_{n'+D}]_{b}\delta_{n'-n'',D'}$. The matrix $K_M(D,D')=([K_M(D,D')]_{a,b})_{1\le a,b\le\bbd}$ is the controlled coefficient matrix for the simple index case. We compute the sum based on the $L^1$ norm of $D,D'\in\Z^{\bbd}$. We get the kernel coefficients as follows:
\begin{itemize}
	\item If $D=D'=0$, we get from \eqref{sim1} and \eqref{sim4}
		\begin{equation*}
			[K_{M}(0,0)]_{b,a}=\frac{3}{2}\sum_{\substack{\gamma^1,\gamma^2,d_1,d_2\\\gamma^{1}_{d_1}=\gamma^2_{d_2}}}[M_{\gamma^1}(d_1)]_{a}[M_{\gamma^2}(d_2)]_{b}.
		\end{equation*}
	\item If $1\le\|D\|_{L^1}\le N$ and $D'=0$, we also get from \eqref{sim1} and \eqref{sim4}
		\begin{equation*}
			[K_{M}(D,0)]_{b,a}=-\sum_{\substack{\gamma^1,\gamma^2\\\gamma^{1}_{\|D\|_{L^1}}=\gamma^2_{\|D\|_{L^1}}=D}}[M_{\gamma^1}(\|D\|_{L^1})]_{a}[M_{\gamma^2}(\|D\|_{L^1})]_{b}-\frac{1}{2}\sum_{\substack{\gamma^1,\gamma^2,d_1,d_2\\\gamma^{1}_{d_1}+\gamma^2_{d_2}=D}}[M_{\gamma^1}(d_1)]_{a}[M_{\gamma^2}(d_2)]_{b}.
		\end{equation*}
	\item If $1\le\|D'\|_{L^1}\le N$ and $D=0$, we obtain the coefficient from \eqref{sim1}, \eqref{sim2}, \eqref{sim3}
		\begin{equation*}
			[K_{M}(0,D')]_{b,a}=-\sum_{\substack{\gamma^1,\gamma^2\\\gamma^{1}_{\|D'\|_{L^1}}=\gamma^2_{\|D'\|_{L^1}}=D'}}[M_{\gamma^1}(\|D'\|_{L^1})]_{a}[M_{\gamma^2}(\|D'\|_{L^1})]_{b}-\frac{1}{2}\sum_{\substack{\gamma^1,\gamma^2,d_1,d_2\\\gamma^{1}_{d_1}+\gamma^2_{d_2}=D'}}[M_{\gamma^1}(d_1)]_{a}[M_{\gamma^2}(d_2)]_{b}.
		\end{equation*}
	\item If $N+1\le\|D\|_{L^1}\le N$ and $D'=0$, the coefficient only involves \eqref{sim4}
		\begin{equation*}
			[K_{M}(D,0)]_{b,a}=-\frac{1}{2}\sum_{\substack{\gamma^1,\gamma^2,d_1,d_2\\\gamma^{1}_{d_1}+\gamma^2_{d_2}=D}}[M_{\gamma^1}(d_1)]_{a}[M_{\gamma^2}(d_2)]_{b}.
		\end{equation*}
	\item If $N+1\le\|D'\|_{L^1}\le N$ and $D=0$, the coefficient only involves \eqref{sim1}
		\begin{equation*}
			[K_{M}(0,D')]_{b,a}=-\frac{1}{2}\sum_{\substack{\gamma^1,\gamma^2,d_1,d_2\\\gamma^{1}_{d_1}+\gamma^2_{d_2}=D'}}[M_{\gamma^1}(d_1)]_{a}[M_{\gamma^2}(d_2)]_{b}.
		\end{equation*}
	\item If $1\le\|D\|_{L^1}\le N$ and $D'=D$, we consider \eqref{sim1}, \eqref{sim2}, \eqref{sim3} and get
		\begin{equation*}
			[K_{M}(D,D)]_{b,a}= -\sum_{\substack{\gamma^1,\gamma^2\\\gamma^{1}_{\|D\|_{L^1}}=\gamma^2_{\|D\|_{L^1}}=D}}[M_{\gamma^1}(\|D\|_{L^1})]_{a}[M_{\gamma^2}(\|D\|_{L^1})]_{b}-\frac{1}{2}\sum_{\substack{\gamma^1,\gamma^2,d_1,d_2\\\gamma^{1}_{d_1}+\gamma^2_{d_2}=D}}[M_{\gamma^1}(d_1)]_{a}[M_{\gamma^2}(d_2)]_{b}.
		\end{equation*}
	\item If $1\le\|D\|_{L^1}\le N/2$ and $D'=-D$, we also consider \eqref{sim1}, \eqref{sim2}, \eqref{sim3} and get
		\begin{equation*}
			\begin{split}
				[K_{M}(D,-D)]_{b,a}=& \sum_{\substack{\gamma^1,\gamma^2\\\gamma^{1}_{\|D\|_{L^1}}=\gamma^2_{\|D\|_{L^1}}=D}}[M_{\gamma^1}(\|D\|_{L^1})]_{a}[M_{\gamma^2}(\|D\|_{L^1})]_{b}\\
				&-\sum_{\substack{\gamma^1,\gamma^2\\\gamma^{1}_{\|D\|_{L^1}}=D,\gamma^2_{2\|D\|_{L^1}}=2D}}[M_{\gamma^1}(\|D\|_{L^1})]_{a}[M_{\gamma^2}(2\|D\|_{L^1})]_{b}\\
				&-\sum_{\substack{\gamma^1,\gamma^2\\\gamma^{1}_{2\|D\|_{L^1}=2D,\gamma^2_{\|D\|_{L^1}}=D}}}[M_{\gamma^1}(2\|D\|_{L^1})]_{a}[M_{\gamma^2}(\|D\|_{L^2})]_{b}.
			\end{split}
		\end{equation*}
	\item If $N/2<\|D\|_{L^1}\le N$ and $D'=- D$, we have
		\begin{equation*}
			[K_{M}(D,-D)]_{b,a}= \sum_{\substack{\gamma^1,\gamma^2\\\gamma^{1}_{\|D\|_{L^1}}=\gamma^2_{\|D\|_{L^1}}=D}}[M_{\gamma^1}(\|D\|_{L^1})]_{a}[M_{\gamma^2}(\|D\|_{L^1})]_{b}.
		\end{equation*}
	\item If $N+1\le\|D\|_{L^1}\le 2N$ and $D'=D$, we consider only \eqref{sim1} and get
		\begin{equation*}
			[K_{M}(D,D)]_{b,a}= -\frac{1}{2}\sum_{\substack{\gamma^1,\gamma^2,d_1,d_2\\\gamma^{1}_{d_1}+\gamma^2_{d_2}=D}}[M_{\gamma^1}(d_1)]_{a}[M_{\gamma^2}(d_2)]_{b}.
		\end{equation*}
	\item If $1\le\|D\|_{L^1}\le N<\|D'\|_{L^1}\le 2N$ and $1\le\|D-D'\|_{L^1}\le N$, only \eqref{sim1} contributes to the coefficient, we compute
		\begin{equation*}
			\begin{split}
				[K_{M}(D,D')]_{b,a}=& \frac{1}{2}\sum_{\substack{\gamma^1,\gamma^2\\\gamma^{1}_{\|D\|_{L^1}}=D,\gamma^2_{\|D'-D\|_{L^1}}=D'-D}}[M_{\gamma^1}(\|D\|_{L^1})]_{a}[M_{\gamma^2}(\|D'-D\|_{L^1})]_{b}\\
				&+\frac{1}{2}\sum_{\substack{\gamma^1,\gamma^2\\\gamma^{1}_{\|D'-D\|_{L^1}}=D'-D,\gamma^2_{\|D\|_{L^1}}=D}}[M_{\gamma^1}(\|D'-D\|_{L^1})]_{a}[M_{\gamma^2}(\|D\|_{L^1})]_{b}.
			\end{split}
		\end{equation*}
	\item If $1\le\|D'\|_{L^1}\le N<\|D\|_{L^1}\le 2N$ and $1\le\|D'-D\|_{L^1}\le N$, from \eqref{sim2} and \eqref{sim3}, we compute
		\begin{equation*}
			\begin{split}
				[K_{M}(D,D')]_{b,a}=& \frac{1}{2}\sum_{\substack{\gamma^1,\gamma^2\\\gamma^{1}_{\|D'\|_{L^1}}=D',\gamma^2_{\|D-D'\|_{L^1}}=D-D'}}[M_{\gamma^1}(\|D'\|_{L^1})]_{a}[M_{\gamma^2}(\|D-D'\|_{L^1})]_{b}\\
				&+\frac{1}{2}\sum_{\substack{\gamma^1,\gamma^2\\\gamma^{1}_{\|D-D'\|_{L^1}}=D-D',\gamma^2_{\|D'\|_{L^1}}=D'}}[M_{\gamma^1}(\|D-D'\|_{L^1})]_{a}[M_{\gamma^2}(\|D'\|_{L^1})]_{b}.
			\end{split}
		\end{equation*}
	\item If $1\le \|D\|_{L^1},\|D'\|_{L^1}\le N$ and $D\ne \pm D'$, we consider \eqref{sim1}, \eqref{sim2}, \eqref{sim3} and get
		\begin{equation*}
			\begin{split}
				[K_{M}(D,D')]_{b,a}=&\frac{1}{2}\sum_{\substack{\gamma^1,\gamma^2\\\gamma^{1}_{\|D\|_{L^1}}=D\\\gamma^2_{\|D'-D\|_{L^1}}=D'-D}}[M_{\gamma^1}(\|D\|_{L^1})]_{a}[M_{\gamma^2}(\|D'-D\|_{L^1})]_{b}\\
				&+\frac{1}{2}\sum_{\substack{\gamma^1,\gamma^2\\\gamma^{1}_{\|D'-D\|_{L^1}}=D'-D\\\gamma^2_{\|D\|_{L^1}}=D}}[M_{\gamma^1}(\|D'-D\|_{L^1})]_{a}[M_{\gamma^2}(\|D\|_{L^1})]_{b}\\
				&+\frac{1}{2}\sum_{\substack{\gamma^1,\gamma^2\\\gamma^{1}_{\|D'\|_{L^1}}=D'\\\gamma^2_{\|D-D'\|_{L^1}}=D-D'}}[M_{\gamma^1}(\|D'\|_{L^1})]_{a}[M_{\gamma^2}(\|D-D'\|_{L^1})]_{b}\\
				&+\frac{1}{2}\sum_{\substack{\gamma^1,\gamma^2\\\gamma^{1}_{\|D-D'\|_{L^1}}=D-D'\\\gamma^2_{\|D'\|_{L^1}}=D'}}[M_{\gamma^1}(\|D-D'\|_{L^1})]_{a}[M_{\gamma^2}(\|D'\|_{L^1})]_{b}\\
				&-\frac{1}{2}\sum_{\substack{\gamma^1,\gamma^2\\\gamma^{1}_{\|D'\|_{L^1}}=D'\\\gamma^2_{\|D\|_{L^1}}=D}}[M_{\gamma^1}(\|D'\|)]_a[M_{\gamma^2}(\|D\|)]_b\\
				&-\frac{1}{2}\sum_{\substack{\gamma^1,\gamma^2\\\gamma^{1}_{\|D\|_{L^1}}=D\\\gamma^2_{\|D'\|_{L^1}}=D'}}[M_{\gamma^1}(\|D\|)]_a[M_{\gamma^2}(\|D'\|)]_b.
			\end{split}
		\end{equation*}
	\item $K_{M}(D,D')=0$ in all other cases.
\end{itemize}
Using all the previous computations, we obtain the coefficients $K_M$ as in Definition \ref{simdefcoeff}.

\subsection{The dual indices case} As in the previous computations, the kernel coefficients come from $\sum_{\gamma}\sum_{n,n',n''}(\lambda^{a,b}_{n,\gamma}\phi_{n'})\cdot(\lambda^{a,b}_{n,\gamma}\phi_{n''}^{*})$ for $a\ne b$ (note that $\lambda^{a,a}_{n,\gamma}=0$).

We have
\begin{align}
		&\sum_{n}(\lambda^{a,b}_{n,\gamma}\phi_{n'})\cdot(\lambda^{a,b}_{n,\gamma}\phi_{n''}^{*})\nonumber\\
		&=  \frac{1}{2}\left( \left( \sum_{d=1}^{N}M^{a,b}_{\gamma}(d)([\alpha_{n+\gamma_d}]_b-[\alpha_n]_{b}) \right)^2+\left( \sum_{d=1}^{N}M^{a,b}_{\gamma}(d)([\alpha_{n+\gamma_d}]_a-[\alpha_n]_{a}) \right)^2 \right)\label{1stdu}\\
		&\times\left( \sum_{d=1}^{N}M_{\gamma}^{a,b}(d)\delta_{n'-n,\gamma_{d}}-\sum_{d=1}^{N}M_{\gamma}^{a,b}(d)\delta_{n,n'} \right)\left( \sum_{d=1}^{N}M_{\gamma}^{a,b}(d)\delta_{n''-n,\gamma_{d}}-\sum_{d=1}^{N}M_{\gamma}^{a,b}(d)\delta_{n,n''} \right).  \label{2nddu}
	\end{align}
Thus, we get the sum
\begin{equation*}
	\sum_{n',n'',D,D'}K^{a,b}(D,D')([\alpha_{n'}]_{a}[\alpha_{n'+D}]_{a}+[\alpha_{n'}]_{b}[\alpha_{n'+D}]_{b})\delta_{n'-n'',D'},
\end{equation*}
where $K^{a,b}$ is computed as follows:
\begin{itemize}
	\item If $D=D'=0$, we first take the coefficients of $[\alpha_n]_a^2$ and $[\alpha_{n+\gamma_d}]_a^2$ in \eqref{1stdu}. Then, we multiply them with the coefficients of pair $\delta_{n,n'}$, $\delta_{n,n''}$ and pairs $\delta_{n'-n,\gamma_d},\delta_{n''-n,\gamma_d}$ having same $d$ in \eqref{2nddu}. We get the coefficient
		\begin{equation*}
			\begin{split}
				K^{a,b}(0,0)=&\frac{1}{2}\sum_{\gamma}\left( \sum_{d=1}^{N}M_{\gamma}^{a,b}(d) \right)^{4}+\sum_{\gamma}\left( \sum_{d=1}^{N}M_{\gamma}^{a,b}(d) \right)^{2}\left( \sum_{d=1}^{N}(M^{a,b}_{\gamma}(d))^2 \right)\\
				&+\frac{1}{2}\sum_{\gamma}\left( \sum_{d=1}^{N}(M^{a,b}_{\gamma}(d))^2 \right)^{2}\\
				=&\frac{1}{2}\sum_{\gamma}\left( \left( \sum_{d=1}^{N}M_{\gamma}^{a,b}(d) \right)^{2}+\left( \sum_{d=1}^{N}(M^{a,b}_{\gamma}(d))^2 \right) \right)^{2}.
			\end{split}
		\end{equation*}
	\item If $1\le \|D\|_{L^1}\le N$ and $D'=0$, we also take the coefficients of pair $\delta_{n,n'}$, $\delta_{n,n''}$ and pairs $\delta_{n'-n,\gamma_d},\delta_{n''-n,\gamma_d}$ having same $d$ in \eqref{2nddu}. Then, we multiply them with the coefficients of $[\alpha_n]_a[\alpha_{n+D}]_a$ in \eqref{1stdu}. We get
		\begin{equation*}
			\begin{split}
				K^{a,b}(D,0)=&-\sum_{\substack{\gamma\\\gamma_{\|D\|_{L^1}}=D}}M^{a,b}_{\gamma}(\|D\|_{L^1})\left(\sum_{d=1}^{N}M^{a,b}_{\gamma}(d)\right)\left( \left( \sum_{d=1}^{N}M_{\gamma}^{a,b}(d) \right)^{2}+\left( \sum_{d=1}^{N}(M^{a,b}_{\gamma}(d))^2 \right) \right)\\
				&+\sum_{\substack{\gamma,d_1,d_2\\\gamma_{d_2}-\gamma_{d_1}=D}}M^{a,b}_{\gamma}(d_1)M^{a,b}_{\gamma}(d_2)\left( \left( \sum_{d=1}^{N}M_{\gamma}^{a,b}(d) \right)^{2}+\left( \sum_{d=1}^{N}(M^{a,b}_{\gamma}(d))^2 \right) \right).
			\end{split}
		\end{equation*}
	\item If $D=0$ and $1\le\|D'\|_{L^1}\le N$, we take the coefficients of $[\alpha_n]_a^2$ and $[\alpha_{n+\gamma_d}]_a^2$ in \eqref{1stdu} (take the factor 1/2 for symmetry between $\gamma$, $-\gamma$). Then, we multiply them with the coefficients of pairs $\delta_{n,n'}$, $\delta_{n''-n,-D'}$, pairs $\delta_{n'-n,D'},\delta_{n,n''}$ and pairs $\delta_{n'-n,\gamma_{d_1}},\delta_{n''-n,\gamma_{d_2}}$ with $\gamma_{d_1}-\gamma_{d_2}=D'$ in \eqref{2nddu}. We obtain
		\begin{equation*}
			\begin{split}
				K^{a,b}(0,D')=&-\frac{1}{2}\sum_{\substack{\gamma\\\gamma_{\|D'\|_{L^1}}=D'}} \left( \left( \sum_{d=1}^{N}M_{\gamma}^{a,b}(d) \right)^{2}+\left( \sum_{d=1}^{N}(M^{a,b}_{\gamma}(d))^2 \right) \right)\left( \sum_{d=1}^{N}M^{a,b}_{\gamma}(d) \right)M^{a,b}_{\gamma}(\|D'\|_{L^1})\\
				&-\frac{1}{2}\sum_{\substack{\gamma\\\gamma_{\|D'\|_{L^1}}=D'}}\left( \left( \sum_{d=1}^{N}M_{-\gamma}^{a,b}(d) \right)^{2}+\left( \sum_{d=1}^{N}(M^{a,b}_{-\gamma}(d))^2 \right) \right)\left( \sum_{d=1}^{N}M^{a,b}_{-\gamma}(d) \right)M^{a,b}_{-\gamma}(\|D'\|_{L^1})\\
				&+\frac{1}{2}\sum_{\substack{\gamma,d_1,d_2\\\gamma_{d_2}-\gamma_{d_1}=D'}}\left( \left( \sum_{d=1}^{N}M_{\gamma}^{a,b}(d) \right)^{2}+\left( \sum_{d=1}^{N}(M^{a,b}_{\gamma}(d))^2 \right) \right)M^{a,b}_{\gamma}(d_1)M^{a,b}_{\gamma}(d_2)\\
				&+\frac{1}{2}\sum_{\substack{\gamma,d_1,d_2\\\gamma_{d_2}-\gamma_{d_1}=D'}}\left( \left( \sum_{d=1}^{N}M_{-\gamma}^{a,b}(d) \right)^{2}+\left( \sum_{d=1}^{N}(M^{a,b}_{-\gamma}(d))^2 \right) \right)M^{a,b}_{\gamma}(d_1)M^{a,b}_{\gamma}(d_2).
			\end{split}
		\end{equation*}
	\item If $1\le\|D\|_{L^1},\|D'\|_{L^1}\le N$ and $D,D'$ are both on a same set $p$ (or $D$ on set $p$ and $D'$ on set $p'$), we get
		\begin{equation*}
			\begin{split}
				K^{a,b}(D,D')=&\sum_{\substack{\gamma\\\gamma_{\|D\|_{L^1}=D}\\\gamma_{\|D'\|_{L^1}}=D'}}\left( \sum_{d=1}^{N}M_{\gamma}^{a,b}(d) \right)^{2}M^{a,b}_{\gamma}(\|D\|_{L^1})M^{a,b}_{\gamma}(\|D'\|_{L^1})\\
				&-\sum_{\substack{\gamma,d_1,d_2\\\gamma_{\|D\|_{L^1}}=D\\\gamma_{d_2}-\gamma_{d_1}=D'}}\left(M^{a,b}_{\gamma}(\|D\|_{L^1})\right)\left(\sum_{d=1}^{N}M^{a,b}_{\gamma}(d)\right)M^{a,b}_{\gamma}(d_1)M^{a,b}_{\gamma}(d_2).
			\end{split}
		\end{equation*}
		We change $\gamma$ to $-\gamma$ for $M^{a,b}_{\gamma}$ in case of $D$ on $p$ and $D'$ on $p'$.
	\item $K^{a,b}(D,D')=0$ in all other cases.
\end{itemize}
Therefore, we get the coefficients $K^{a,b}$ as in Definition \ref{dualdefcoeff}.

\section{Proofs of the main results}
Since the extension to the multi-dimensional case is straightforward, to make the following proofs more readable we restrict ourselves, hereafter, to the 1-dimensional case.

\subsection{Proof of Lemma \ref{quotientbound}}
Let us first establish that
\begin{equation}\label{sum01D}
	\sum_{|d|\le 2N} K_{M_N}(d,d')=0\qquad\forall d'\in \N,|d'|\le2N.
	\end{equation}

	When $d'=0$, we have
	\begin{align}
			\sum_{|d|\le 2N} K_{M_N}(d,0)&= \frac{3}{2}\sum_{1\le |d|\le N}M_N(d)^2-\sum_{1\le |d|\le N}M_N(d)^2-\frac{1}{2}\sum_{1\le |d|\le 2N}\sum_{\substack{d_1+d_2=d\\1\le|d_1|,|d_2|\le N}}M_N(d_1)M_N(d_2)\nonumber\\
			&= -\frac{1}{2}\sum_{1\le |d|\le 2N}\sum_{\substack{d_1+d_2=d\\1\le|d_1|,|d_2|\le N\\d_1\ne d_2}}M_N(d_1)M_N(d_2).\label{d'=0sum}
		\end{align}
		The sum \eqref{d'=0sum} is $0$ because if $(d_1,d_2)$ is in the sum then so is $(-d_1,d_2)$. Hence, \eqref{sum01D} holds when $d'=0$.

	Now, we consider $1\le d'\le N$ (the case $-N\le d'\le -1$ is similar).
	By using the change of variable $d\to -d$, we sum all coefficients of the case $1\le|d|\le N$ into $\sum_{1\le |d'-d|,|d|\le N}[M_N(d)M_N(d'-d)]$. We get
	\begin{equation*}
		\begin{split}
			\sum_{|d|\le 2N} K_{M_N}(d,d')&= -M_N(d')^2 -\sum_{d_1+d_2=d',1\le|d_1|,|d_2|\le N}M_N(d_1)M_N(d_2)-2M_N(d')M_N(2d')\\
			&+\sum_{d=N+1}^{N+d'}M_N(d')M_N(d-d')+\sum_{1\le |d'-d|,|d|\le N} M_N(d)M_N(d'-d)\\
			&-\sum_{1\le |d'-d|,|d|\le N} M_N(d')M_N(d'-d).
		\end{split}
	\end{equation*}

	If $d'> N/2$, for each $d_1\in[N+1,N+d']\setminus\{2d'\}$, we can choose $d_2=2d'-d_1\in[d'-N,2d'-N-1]\setminus\{0\}$. In case $d_1=2d'$, it gives $M_N(d')^2$ which cancels $-M_N(d')^2$. The remainder is
	\begin{equation*}
		-\sum_{\substack{1\le |d'-d|,|d|\le N\\|2d'-d|\le N}}M_N(d')M_N(d'-d).
	\end{equation*}
	Note that if $d$ is in the sum then $2d'-d$ is also in the sum. Since $d'-d$ and $d'-(2d'-d)$ are opposite one to another the sum is $0$. 

	If $d'\le N/2$, for each $d_1\in[N+1,N+d']$, we can choose $d_2=2d'-d_1\in[d'-N,2d'-N-1]\subset[-N,-1]$. Those terms cancel each other. The remainder is 
	\begin{equation}
\label{remain<N/2}
		-M_N(d')^2-2M_N(d')M_N(2d')-\sum_{\substack{1\le |d'-d|,|d|\le N|2d'-d|\le N}}M_N(d')M_N(d'-d).
	\end{equation}
	$M_N(d')M_N(2d')$ is $0$ thanks to Assumption \ref{N/2con}. The case $d=2d'$ makes a term that cancels $-M_N(d')^2$. For $d\ne 2d'$, $d$ and $2d'-d$ make two opposite terms. This leads to \eqref{remain<N/2}$\,=0$.

Hence, \eqref{sum01D} holds when $|d'|\le N$.

Finally, we consider $N+1\le d'\le 2N$ (the case $-2N\le d'\le -N-1$ is similar). We have
	\begin{equation*}
		\sum_{|d|\le 2N}K_{M_N}(d,d')=-\sum_{\substack{d_1+d_2=d'\\1\le|d_i|\le N}}M_N(d_1)M_N(d_2)+\sum_{d=d'-N}^{N}M(d)M_N(d'-d)=0.
	\end{equation*}
	This concludes the proof of \eqref{sum01D}.

	Using the fact that $K_{M_N}(0,d')=-\sum_{1\le |d|\le 2N}K_{M_N}(d,d')$, we can rewrite $\hat{K}_{M_N}(k,k')$ as
\begin{equation*}
	\sum_{|d'|\le 2N}\cos(2\pi d'k')\left( \sum_{1\le|d|\le 2N}K_{M_N}(d,d')(\cos(2\pi dk)-1) \right).
\end{equation*}
Since $|\cos(2\pi nk)-1|\lesssim |k|$ we get $\hat{K}_{M_N}(k,k')\lesssim |k|$. We also see that $|k|\lesssim \omega(k)$. In the pinning case, we already have $|k|\lesssim 1\lesssim \omega(k)$. If $|k|$ is not close to $0$ in the no-pinning case, we still has $|k|\lesssim 1\lesssim \omega(k)$. We consider the no-pinning case and suppose that $|k|$ is close to $0$, we have
\begin{equation*}
	\hat{\sigma}(k)=\hat{\sigma}(0)+\hat{\sigma}'(0)k+\frac{1}{2}\hat{\sigma}''(k')k^2=\frac{1}{2}\hat{\sigma}''(k')k^2,
\end{equation*}
where $0\le|k'|<|k|$. Since $\hat{\sigma}''(0)>0$, when $|k|$ is close to $0$ we also have $\hat{\sigma}''(k')\approx 1$. Thus, $\omega(k)=\sqrt{\hat{\sigma}(k)}\approx |k|$ when $k$ is close to $0$.
Therefore,
\begin{equation*}
	\sup\left|\frac{\hat{K}_{M_N}(k,k')}{\omega(k)}\right|<\infty.
\end{equation*}

\subsection{Proof of Theorem \ref{Enerlem}}
\subsubsection{The energy transport equation \eqref{energytransport}}
		Let $k\mapsto S(k)$ be a bounded real-valued function on $\T$, we have
		\begin{equation*}
			\frac{d}{dt}\langle S,E^{\varepsilon}(t)\rangle= \frac{\varepsilon}{2}\int_{\T}\varepsilon^{-1}\E_{\varepsilon}\left[ O|\hat{\phi}(k,t/\varepsilon)|^2 \right]S(k)dk.
		\end{equation*}
		We use the techniques developed in Section \ref{comker} and compute
		\begin{equation*}
			\begin{split}
				\sum_{n,n',n''\in\Z}\left( \lambda^{M_N}_n\phi_{n'} \right)\left( \lambda^{M_N}_n\phi_{n''}^{*} \right)\tilde{S}^{*}(n''-n')&= \sum_{n',d,d'}K_{M_N}(d,d')\alpha_{n'}\alpha_{n'+d}\tilde{S}^{*}(-d')\\
				&= \int_{\T}|\hat{\alpha}(k)|^2\sum_{d,d'}e^{-2\pi i kd}K_{M_N}(d,d')\tilde{S}^{*}(-d')dk\\
				&= \int_{\T^2}|\hat{\alpha}(k)|^2\hat{K}_{M_N}(k,k')S(k')dkdk'.
			\end{split}
		\end{equation*}
Thus, we get
\begin{equation*}
	\int_{\T}O^{M_N}|\hat{\phi}(k)|^2S(k)dk=\int_{\T}\left( |\hat{\phi}(k)|^2-\frac{1}{2}\left( \hat{\phi}(k)\hat{\phi}(-k)+\hat{\phi}^{*}(k)\hat{\phi}^{*}(-k) \right) \right)O^{M_N}_{col}S(k)dk.
\end{equation*}
Therefore, the time derivative is
\begin{equation*}
	\begin{split}
		\frac{d}{dt}\langle S,E^{\varepsilon}(t)\rangle&= \langle O^{M_N}_{col}S,E^{\varepsilon}(t)\rangle-\frac{\varepsilon}{4}\int_{\T}\left( \hat{\phi}(k,t/\varepsilon)\hat{\phi}(-k,t/\varepsilon)+\hat{\phi}^{*}(k,t/\varepsilon)\hat{\phi}^{*}(-k,t/\varepsilon) \right)O^{M_N}_{col}S(k)dk.
	\end{split}
\end{equation*}
It is sufficient to check that
\begin{equation}
\label{limprod=0}
	\lim_{\varepsilon\to}\int_{0}^{t}\frac{\varepsilon}{2}\int_{\T}\hat{\phi}(k,s/\varepsilon)\hat{\phi}(-k,s/\varepsilon)O^{M_N}_{col}S(k)dkds=0, \forall t\in[0,T].
\end{equation}
The complex conjugate counterpart is proved similarly. To prove \eqref{limprod=0}, we compute the derivative
\begin{equation*}
	\frac{d}{dt}\frac{\varepsilon}{2}\int_{\T}\hat{\phi}(k,t/\varepsilon)\hat{\phi}(-k,t/\varepsilon)S(k)dk=\frac{\varepsilon}{2}\int_{\T}\varepsilon^{-1}O\left( \hat{\phi}(k,t/\varepsilon)\hat{\phi}(-k,t/\varepsilon) \right)S(k)dk.
\end{equation*}
The Hamiltonian operator gives
\begin{equation*}
	O^H\left( \hat{\phi}(k)\hat{\phi}(-k) \right)=-2i\omega(k)\hat{\phi}(k)\hat{\phi(-k)}.
\end{equation*}
We use the techniques developed in Section \ref{comker} again and see that
\begin{equation*}
	\begin{split}
		\int_{\T}O^{M_N}\left( \hat{\phi}(k)\hat{\phi}(-k) \right)S(k)dk&=  -\frac{d}{dt}\langle S,E^{\varepsilon}(t)\rangle\\
		&+ \frac{\varepsilon}{4}\int_{\T^2}\hat{\phi}(k,t/\varepsilon)\hat{\phi}(-k,t/\varepsilon)\hat{K}_{M_N}(k,k')S(k)dkdk'\\
		&- \frac{\varepsilon}{4}\int_{\T^2}\hat{\phi}^{*}(k,t/\varepsilon)\hat{\phi}^{*}(-k,t/\varepsilon)\hat{K}_{M_N}(k,k')S(k)dkdk'.
	\end{split}
\end{equation*}
Hence, we get
\begin{equation*}
	\begin{split}
		\frac{d}{dt}\frac{\varepsilon}{2}\int_{\T}\hat{\phi}(k,t/\varepsilon)\hat{\phi}(-k,t/\varepsilon)S(k)dk&= -\frac{2i}{\varepsilon}\frac{\varepsilon}{2}\int_{\T}\hat{\phi}(k,t/\varepsilon)\hat{\phi}(-k,t/\varepsilon)\omega(k)S(k)dk-\frac{d}{dt}\langle S,E^{\varepsilon}(t)\rangle\\
		&+ \frac{\varepsilon}{4}\int_{\T^2}\hat{\phi}(k,t/\varepsilon)\hat{\phi}(-k,t/\varepsilon)\hat{K}_{M_N}(k,k')S(k)dkdk'\\
		&- \frac{\varepsilon}{4}\int_{\T^2}\hat{\phi}^{*}(k,t/\varepsilon)\hat{\phi}^{*}(-k,t/\varepsilon)\hat{K}_{M_N}(k,k')S(k)dkdk'.
	\end{split}
\end{equation*}
Assumption \ref{boundenergy} ensures that $\frac{\varepsilon}{2}\int_{\T}\hat{\phi}(k,t/\varepsilon)\hat{\phi}(-k,t/\varepsilon)S(k)dk$, $\frac{\varepsilon}{2}\int_{\T}\hat{\phi}^{*}(k,t/\varepsilon)\hat{\phi}^{*}(-k,t/\varepsilon)S(k)dk$, $\langle S,E^{\varepsilon}(t)\rangle$ are all bounded. By integrating with respect to $t$, we obtain
\begin{equation*}
	\lim_{\varepsilon\to0}\int_{0}^{t}\frac{\varepsilon}{2}\int_{\T}\hat{\phi}(k,s/\varepsilon)\hat{\phi}(-k,s/\varepsilon)\omega(k)S(k)dkds=0.
\end{equation*}
Lemma \ref{quotientbound} implies that $|\hat{K}_{M_N}(k,k')|\lesssim \omega(k)$ and therefore we can replace $S$ by $\frac{O^{M_N}_{col}S(k)}{\omega(k)}$. The equation \eqref{energytransport} is proved.

	\subsubsection{Proof of \eqref{nopinex}}
		From the proof of \eqref{energytransport}, we have
		\begin{equation*}
			\begin{split}
				\left|\frac{d}{dt}\langle 1_{(-R,R)},E^{\varepsilon}(t)\rangle\right|&\le \left|\langle O^{M_N}_{col}1_{(-R,R)},E^{\varepsilon}(t)\rangle\right|\\
				&+\left|\frac{\varepsilon}{4}\int_{\T}\left( \hat{\phi}(k,t/\varepsilon)\hat{\phi}(-k,t/\varepsilon)+\hat{\phi}^{*}(k,t/\varepsilon)\hat{\phi}^{*}(-k,t/\varepsilon) \right)O^{M_N}_{col}1_{(-R,R)}dk\right|.
			\end{split}
		\end{equation*}
		Thus, by the fact that $\hat{K}_{M_N}$ is bounded, we get $O^{M_N}_{col}1_{(-R,R)}\lesssim R+1_{(-R,R)}$. We have the estimate
		\begin{equation*}
			\left|\frac{d}{dt}\langle 1_{(-R,R)},E^{\varepsilon}(t)\rangle\right|\lesssim R+\langle 1_{(-R,R)},E^{\varepsilon}(t)\rangle.
		\end{equation*}
		By Gronwall's inequality, we get
		\begin{equation*}
			\langle 1_{[-R,R]},E^{\varepsilon}(t)\rangle\lesssim \left(R+\langle 1_{(-R,R)},E^{\varepsilon}(0)\rangle\right)e^{C_8 t}
		\end{equation*}
		for a constant $C_8>0$ and each $t\in[0,T]$. By Assumption \ref{nopin}, the conclusion follows.

	\subsection{Proof of Theorem \ref{main}}
		We recall the definition of the Wigner distribution in the Fourier space
		\begin{equation*}
			\langle S,W^{\varepsilon}(t)\rangle :=\frac{\varepsilon}{2}\int_{\T_{2/\varepsilon}\times\T}\hat{\phi}(k+\varepsilon\xi/2)\hat{\phi}^{*}(k-\varepsilon\xi/2)\hat{S}^{*}(\xi,k),
		\end{equation*}
		where $\T_{2/\varepsilon}$ is the torus of size $2/\varepsilon$ and $\hat{S}(\xi,k)=\int_{\R}e^{-2\pi i x\xi}S(x,k)$.
		Using the definition, we compute
		\begin{equation*}
			\begin{split}
				\frac{d}{dt}\langle S,W^{\varepsilon}(t)\rangle&= \frac{\varepsilon}{2}\int_{\T_{2/\varepsilon}\times\T}\varepsilon^{-1}\E_{\varepsilon}\left[ O^H(\hat{\phi}(k+\varepsilon\xi/2,t/\varepsilon)\hat{\phi}^{*}(k-\varepsilon\xi/2,t/\varepsilon)) \right]\hat{S}^{*}(\xi,k)d\xi dk\\
				&+\frac{\varepsilon}{2}\sum_{n',n''}\E_{\varepsilon}\left[ O^{M_N}(\phi_{n'}(t/\varepsilon)\phi_{n''}^{*}(t/\varepsilon)) \right]\tilde{S}^{*}\left( \frac{\varepsilon(n'+n'')}{2},n'-n'' \right).
			\end{split}
		\end{equation*}
		First, we deal with the Hamiltonian operator. By \eqref{O^Hiden}, we have
		\begin{equation}
\label{OHcom}
			O^H(\hat{\phi}(k+\varepsilon\xi/2)\hat{\phi}^{*}(k-\varepsilon\xi/2))=-i(\omega(k+\varepsilon\xi/2)-\omega(k-\varepsilon\xi/2))\hat{\phi}(k+\varepsilon\xi/2)\hat{\phi}^{*}(k-\varepsilon\xi/2).
		\end{equation}
		We see that $\varepsilon^{-1}(\omega(k+\varepsilon\xi/2)-\omega(k-\varepsilon\xi/2))$ can be estimated by $\omega'(k)\xi$. We prove that
		\begin{equation}
			\begin{split}
						\frac{\varepsilon}{2}\int_{\T_{2/\varepsilon}\times\T}\E_{\varepsilon}&\left[ \hat{\phi}(k+\varepsilon\xi/2,t/\varepsilon)\hat{\phi}^{*}(k-\varepsilon\xi/2,t/\varepsilon) \right]\\
								&\times(\varepsilon(\omega(k+\varepsilon\xi/2)-\omega(k-\varepsilon\xi/2))-\omega'(k)\xi)\hat{S}^{*}(\xi,k)d\xi dk
							\end{split}
			\label{estomega'}
		\end{equation}
		converges to $0$ when $\varepsilon\to0$. 
		We consider the no-pinning case and we split the domain of integration into three domains: $|k|\ge R>\varepsilon|\xi|$, $R\le\varepsilon|\xi|$ and $|k|<R$ for sufficiently small numbers $\varepsilon,R$.

		In the case $|k|\ge R>\varepsilon|\xi|$, $\omega$ is smooth on the interval between $k\pm \varepsilon\xi/2$ (consider it on $(0,1)$ or $(-1,0)$ as $k$ may be close to $\pm 1/2$). Thus, we get
		\begin{equation}
			\label{o''est}
			|\varepsilon^{-1}(\omega(k+\varepsilon\xi/2)-\omega(k-\varepsilon\xi/2))-\omega'(k)\xi|=\frac{1}{2}|\omega''(k')|\varepsilon\xi^2
		\end{equation}
		where $k'$ is between $k\pm \varepsilon\xi/2$. We have an estimate on the second-order derivative
		\begin{equation*}
			\omega''(k')=\frac{\hat{\sigma}''(k')\omega(k')-\hat{\sigma}'(k')\omega'(k')}{2\hat{\sigma}(k')}=\frac{\hat{\sigma}''(k')\omega(k')-(\hat{\sigma}'(k'))^2/(2\omega(k'))}{2\hat{\sigma}(k')}.
		\end{equation*}
		We see that $|\omega''(k')|\approx|k'|^{-1}\approx |k|^{-1}\lesssim R^{-1}$ because $\hat{\sigma}''(k')\approx 1$, $\omega(k')\approx|k'|$, $\hat{\sigma}'(k')\approx |k'|$ and $\hat{\sigma}(k')\approx |k'|^2$. The integral \eqref{estomega'} on this domain is bounded by
		\begin{equation*}
			\varepsilon\int_{|k|>R}\E_{\varepsilon}\left[ |\hat{\phi}(k)|^2 \right]dk\sup_{k}\left|\int_{-R\varepsilon^{-1}}^{R\varepsilon^{-1}}\hat{S}^{*}(\xi,k)R^{-1}\varepsilon\xi^{2}d\xi\right|.
		\end{equation*}
		As $S$ is in the Schwartz space, we bound $\int_{\R}|\hat{S}(\xi,k)||\xi|^2d\xi$ by a constant. Hence, the integral \eqref{estomega'} on this domain tends to $0$ as $\varepsilon\to0$ for each fixed $R$.

		Considering the case $R\le\varepsilon|\xi|$, we use the estimate
		\begin{equation*}
		|\varepsilon^{-1}(\omega(k+\varepsilon\xi/2)-\omega(k-\varepsilon\xi/2))-\omega'(k)\xi)|\lesssim R^{-1}|\xi|.
		\end{equation*}
		The integral \eqref{estomega'} on this domain is bounded by
	\begin{equation*}
		\varepsilon\int_{\T}\E_{\varepsilon}\left[ |\hat{\phi}(k)|^2 \right]dk\sup_{k}\int_{R\le\varepsilon|\xi|\le2}\hat{S}^{*}(\xi,k)R^{-1}|\xi|d\xi.
		\end{equation*}
		We use $|\xi|^{-1}\le R^{-1}\varepsilon$ and again we bound $\int_{\R}|\hat{S}(\xi,k)||\xi|^2d\xi$. The integral \eqref{estomega'} also tends to $0$ as $\varepsilon\to0$ on this domain.

		For the case $|k|<R$, we use \eqref{nopinex}. We will need to bound
		\begin{equation}\label{k<Rbound}
			\int_{\T_{2/\varepsilon}}\hat{S}^{*}(\xi,k)(\varepsilon^{-1}(\omega(k+\varepsilon\xi/2)-\omega(k-\varepsilon\xi/2))-\omega'(k)\xi)d\xi.
		\end{equation}
		When $|k|>\varepsilon|\xi|$, we use \eqref{o''est} with the estimate $|\omega''(k')|\varepsilon|\xi|\lesssim 1$ and use a bound for $\int_{\R}|\hat{S}||\xi|d\xi$. When $|k|\le \varepsilon|\xi|$ we have
		\begin{equation*}
		|\varepsilon^{-1}(\omega(k+\varepsilon\xi/2)-\omega(k-\varepsilon\xi/2))-\omega'(k)\xi)|\le \varepsilon^{-1}(\varepsilon|\xi|/2+|k|+\varepsilon|\xi|/2-|k|)+|\omega'(k)||\xi|\lesssim |\xi|.
		\end{equation*}
		We again use a bound for $\int_{\R}|\hat{S}||\xi|d\xi$. Thus, \eqref{k<Rbound} is bounded.

		In the pinning case, $\omega$ is smooth, we use \eqref{o''est} and the fact that $\omega''$ is bounded and we easily get the limit $0$ for \eqref{estomega'} when $\varepsilon\to0$.
Therefore, the Hamiltonian operator gives
		\begin{equation*}
			\frac{\varepsilon}{2}\int_{\T_{2/\varepsilon}\times\T}\E_{\varepsilon}\left[ \hat{\phi}(k+\varepsilon\xi/2,t/\varepsilon)\hat{\phi}^{*}(k-\varepsilon\xi/2,t/\varepsilon) \right](-i\xi)\omega'(k)\hat{S}^{*}(\xi,k)d\xi dk.
		\end{equation*}
		Since $\hat{\frac{\partial S}{\partial x}}(\xi,k)=(-2\pi i\xi)\hat{S}(\xi,k)$, the term \eqref{OHcom} at the limit $\varepsilon\to0$ is
		\begin{equation*}
			\frac{1}{2\pi}\left\langle \omega'(k)\frac{\partial S}{\partial x},W^{\varepsilon}(t)\right\rangle.
		\end{equation*}

		The perturbation is computed using the techniques developed in Section \ref{comker}. We write
		\begin{equation*}
			\begin{split}
				\sum_{n',n''}&\E_{\varepsilon}\left[ \lambda^{M_N}_n\phi_{n'}\lambda^{M_N}_n\phi_{n''}^{*} \right]S^{*}\left( \frac{\varepsilon(n'+n'')}{2},n'-n'' \right)\\
				&= \sum_{n',d,d'}\E_{\varepsilon}\left[ \alpha_{n'}\alpha_{n'+d} \right]K_{M_N}(d,d')\tilde{S}^{*}\left( \frac{\varepsilon(2n'+d')}{2},d' \right)\\
				&= \sum_{n',n'',d'}\E_{\varepsilon}\left[ \alpha_{n'}\alpha_{n''} \right]K_{M_N}(n''-n',d')\tilde{S}^{*}\left( \frac{\varepsilon(n'+n'')}{2},d' \right)+a(\varepsilon)\\
				&= \frac{1}{2}\sum_{n',n''}\E_{\varepsilon}\left[ \phi_{n'}\phi_{n''}^{*}+\phi_{n'}^{*}\phi_{n''} \right]\int_{\T^2}e^{2\pi ik(n''-n')}\hat{K}_{M_N}(k,k')S^{*}\left( \frac{\varepsilon(n'+n'')}{2},k' \right)dkdk'\\
				&-\frac{1}{2}\sum_{n',n''}\E_{\varepsilon}\left[ \phi_{n'}\phi_{n''}+\phi_{n'}^{*}\phi_{n''}^{*} \right]\int_{\T^2}e^{2\pi ik(n''-n')}\hat{K}_{M_N}(k,k')S^{*}\left( \frac{\varepsilon(n'+n'')}{2},k' \right)dkdk'+a(\varepsilon)
			\end{split}
		\end{equation*}
		where $a(\varepsilon)$ satisfies $\varepsilon a(\varepsilon)\to0$ when $\varepsilon\to0$. We change $\phi_{n'}^{*}\phi_{n''}$ into $\phi_{n'}\phi_{n''}^{*}$. This is possible because $\hat{K}_{M_N}(k,k')$ is even with respect to $k$. We also get that
		\begin{equation*}
			\begin{split}
				\sum_{n',n''}&\E_{\varepsilon}\left[ (O^{M_N}\phi_{n'})\phi_{n''}^{*}+(O^{M_N}\phi_{n''}^{*})\phi_{n'} \right]\tilde{S}^{*}\left( \frac{\varepsilon(n'+n'')}{2},n'-n'' \right)\\
				&= \sum_{n',n''}\E_{\varepsilon}\left[ (-\sum_{d}K_{M_N}(-d,0)(\phi_{n'+d}-\phi_{n'+d}^{*}))\phi_{n''}^{*}+(\sum_{d}K_{M_N}(d,0)(\phi_{n''-d}-\phi_{n''-d}^{*}))\phi_{n'} \right)\\
				&\hspace{20em}\times\tilde{S}^{*}\left( \frac{\varepsilon(n'+n'')}{2},n'-n'' \right)\\
				&= -2\sum_{n',n''}\E_{\varepsilon}\left[ \phi_{n'}\phi_{n''}^{*}-\frac{1}{2}(\phi_{n'}\phi_{n''}+\phi_{n'}^{*}\phi_{n''}^{*}) \right]\sum_{d}K_{M_N}(d,0)\tilde{S}^{*}\left( \frac{\varepsilon(n'+n'')}{2},n'-n''+d \right)+a(\varepsilon)\\
				&= -2\sum_{n',n''}\E_{\varepsilon}\left[ \phi_{n'}\phi_{n''}^{*}-\frac{1}{2}(\phi_{n'}\phi_{n''}+\phi_{n'}^{*}\phi_{n''}^{*}) \right]\int_{\T^2}e^{2\pi ik(n''-n')}\hat{K}(k,k')S^{*}\left( \frac{\varepsilon(n'+n'')}{2},k \right)dkdk'+a(\varepsilon).
			\end{split}
		\end{equation*}
		Therefore, the perturbation gives
		\begin{equation*}
			\frac{\varepsilon}{2}\sum_{n',n''}\E_{\varepsilon}\left[ \phi_{n'}\phi_{n''}^{*}-\frac{1}{2}(\phi_{n'}\phi_{n''}+\phi_{n'}^{*}\phi_{n''}^{*}) \right]\int_{\T}e^{2\pi ik(n''-n')}O^{M_N}_{col}S^{*}\left( \frac{\varepsilon(n'+n'')}{2},k \right)dk+\varepsilon a(\varepsilon).
		\end{equation*}

		To get the conclusion, we show that
		\begin{equation*}
			\lim_{\varepsilon\to0}\frac{\varepsilon}{2}\int_{\T}\E_{\varepsilon}\left[ \phi_{n}(t/\varepsilon)\phi_{n''}(t/\varepsilon) \right]\int_{\T}e^{2\pi ik(n''-n')}O^{M_N}_{col}S^{*}\left( \frac{\varepsilon(n'+n'')}{2},k \right)dk=0.
		\end{equation*}
		We have
		\begin{equation*}
			\begin{split}
			\frac{d}{dt}\frac{\varepsilon}{2}\int_{\T}&\E_{\varepsilon}\left[ \phi_{n}(t/\varepsilon)\phi_{n''}(t/\varepsilon) \right]\tilde{S}^{*}\left( \frac{\varepsilon(n'+n'')}{2},n'-n'' \right)\\
				=&\ \frac{\varepsilon}{2}\int_{\T_{2/\varepsilon}\times\T}\varepsilon^{-1}\E_{\varepsilon}\left[ O^H(\hat{\phi}(k+\varepsilon\xi/2,t/\varepsilon))\hat{\phi}(k-\varepsilon\xi/2,t/\varepsilon) \right]\hat{S}^{*}(\xi,k)d\xi dk\\
				&+\frac{\varepsilon}{2}\sum_{n',n''}\E_{\varepsilon}\left[ O^{M_N}(\phi_{n'}(t/\varepsilon)\phi_{n''}(t/\varepsilon)) \right]\tilde{S}^{*}\left( \frac{\varepsilon(n'+n'')}{2},n'-n'' \right)\\
				=&\ \frac{-i\varepsilon}{2}\int_{\T_{2/\varepsilon}\times\T}\E_{\varepsilon}\left[ \hat{\phi}(k+\varepsilon\xi/2,t/\varepsilon)\hat{\phi}(k-\varepsilon\xi/2,t/\varepsilon) \right]\varepsilon^{-1}(\omega(k+\varepsilon\xi/2)+\omega(k-\varepsilon\xi/2))\hat{S}^{*}(\xi,k)d\xi dk\\
				&-\frac{1}{2}\left( \frac{d}{dt}\langle S,W^{\varepsilon}(t)\rangle+\frac{d}{dt}\langle S,W^{\varepsilon}(t)^{*}\rangle \right)+\frac{1}{2\pi}\left( \langle\omega'\frac{\partial S}{\partial x},W^{\varepsilon}(t)\rangle+\langle\omega'\frac{\partial S}{\partial x},W^{\varepsilon}(t)^{*}\rangle \right)\\
				&+\frac{\varepsilon}{4}\sum_{n',n''}\E_{\varepsilon}\left[ \phi_{n'}\phi_{n''}+\phi_{n'}^{*}\phi_{n''}^{*} \right]\int_{\T^2}e^{2\pi ik(n''-n')}\hat{K}_{M_N}(k,k')S^{*}\left( \frac{\varepsilon(n'+n'')}{2},k' \right)dkdk'\\
				&-\frac{\varepsilon}{4}\sum_{n',n''}\E_{\varepsilon}\left[ \phi_{n'}\phi_{n''}-\phi_{n'}^{*}\phi_{n''}^{*} \right]\int_{\T^2}e^{2\pi ik(n''-n')}\hat{K}_{M_N}(k,k')S^{*}\left( \frac{\varepsilon(n'+n'')}{2},k \right)dkdk'.
			\end{split}
		\end{equation*}
			We get that
			\begin{equation*}
				\int_{0}^{t}\frac{\varepsilon}{2}\int_{\T_{2/\varepsilon}\times\T}\E_{\varepsilon}\left[ \hat{\phi}(k+\varepsilon\xi/2,s/\varepsilon)\hat{\phi}(k-\varepsilon\xi/2,s/\varepsilon) \right](\omega(k+\varepsilon\xi/2)+\omega(k-\varepsilon\xi/2))\hat{S}^{*}(\xi,k)d\xi dkds
			\end{equation*}
			tends to $0$ when $\varepsilon\to0$ for each $t\in[0,T]$. We replace $\omega(k+\varepsilon\xi/2)+\omega(k-\varepsilon\xi/2)$ by $2\omega(k)$. This is feasible because if $|k|>\varepsilon|\xi|$ then
			\begin{equation*}
				\omega(k+\varepsilon\xi/2)+\omega(k-\varepsilon\xi/2)-2\omega(k)=\frac{1}{2}(\omega''(k_+)+\omega''(k_-))\varepsilon^{2}\xi^2.
			\end{equation*}
			In the no-pinning case $\omega''(k_{\pm})\varepsilon|\xi|\lesssim 1$ and in the pinning case $\omega''(k_{\pm})\lesssim 1$. Thus, in both cases, there is always at least one free $\varepsilon$. If $|k|\le\varepsilon|\xi|$ then
\begin{equation*}
	|\omega(k+\varepsilon\xi/2)+\omega(k-\varepsilon\xi/2)-2\omega(k)|\lesssim\varepsilon|\xi|/2+|k|+\varepsilon|\xi|/2-|k|+2|k|\lesssim \varepsilon|\xi|.
			\end{equation*}
			We also have a free $\varepsilon$. We replace $\omega(k+\varepsilon\xi/2)+\omega(k-\varepsilon\xi/2)$ by $2\omega(k)$ and get that
			\begin{equation*}
				\int_{0}^{t}\frac{\varepsilon}{2}\int_{\T_{2/\varepsilon}\times\T}\E_{\varepsilon}\left[ \hat{\phi}(k+\varepsilon\xi/2,s/\varepsilon)\hat{\phi}(k-\varepsilon\xi/2,s/\varepsilon) \right]2\omega(k)\hat{S}^{*}(\xi,k)d\xi dkds
			\end{equation*}
			tends to $0$ when $\varepsilon\to0$ for $t\in[0,T]$. Replacing $S$ by $\frac{O^{M_N}_{col}S(k)}{\omega(k)}$ we get the desired result.

	\subsection{Proof of Proposition \ref{propd1}}
	In the 1-dimensional case, given any $N$, we have
	\begin{equation*}
		\begin{split}
			\hat{K}_{M_N}(k,k')=& 4\sum_{1\le |d|,|d'|\le 2N}K_{M_N}(d,d')\sin^2(\pi dk)\sin^2(\pi d'k')\\
			=& 16\sum_{1\le d, d'\le 2N}K_{M_N}(d,d')\sin^2(\pi k)\sin^2(\pi k')\left( \sum_{d_1=0}^{d-1}U_{2d_1}(\cos(\pi k)) \right)\left( \sum_{d_2=0}^{d'-1}U_{2d_2}(\cos(\pi k')) \right)
		\end{split}
	\end{equation*}
	where $U_d(x)$ is the Chebyshev polynomial defined by
	\begin{equation*}
		\sin((d+1)\theta)=\sin(\theta)U_{d}(\cos(\theta)),\quad d\ge0.
	\end{equation*}
	Noting that
	\begin{equation*}
		\begin{split}
			\sin^2(\theta)\left(\sum_{d_1=0}^{d-1}U_{2d_1}(\cos(\theta))\right)=& \sum_{d_1=0}^{d-1}\sin(\theta)\sin((2d_1-1)\theta)\\
			=& \sum_{d_1=0}^{d-1}\frac{\cos(2d_1\theta)-\cos( 2(d_1+1)\theta)}{2}\\
			=& \frac{1-\cos(2d\theta)}{2}=\sin^2(d\theta)
		\end{split}
	\end{equation*}
	we see that the kernel $\hat{K}_{M_N}$ is of the form
\begin{equation*}
	16\sin^2(\pi k)\sin^2(\pi k')\sum_{d_1=0}^{2N-1}\sum_{d_2=0}^{2N-1}L_{M_N}(d_1,d_2)U_{2d_1}(\cos(\pi k))U_{2d_2}(\cos(\pi k')),
\end{equation*}
where $L_{M_N}$ is defined by
\begin{equation*}
	L_{M_N}(d_1,d_2)=\sum_{d=d_1+1}^{2N}\sum_{d'=d_2+1}^{2N}K_{M_N}(d,d').
	\end{equation*}
	According to \eqref{00} to \eqref{sum<=N} and Assumption \ref{N/2con} (pick $M_N(d)=0,d\le N/2$ for simplicity), $L_{M_N}(d_1,d_2)$ can be rewritten as the sum
	\begin{equation*}
		\sum_{N/2<d\le d'\le 2N} c_{d,d'}(d_1,d_2)M_N(d)M_N(d').
	\end{equation*}
	We define the polymial $P_{i,j}(x,y)$ by
	\begin{equation*}
		P_{i,j}(x,y)=\sum_{d_1,d_2}c_{i+\lfloor N/2\rfloor,j+\lfloor N/2\rfloor}(d_1,d_2)U_{2d_1}(x)U_{2d_2}(y).
	\end{equation*}
	Therefore, for
	\begin{equation*}
		\mathcal{K}(k,k')=\sin^2(\pi k)\sin^2(\pi k')v(\cos(\pi k),\cos(\pi,k'))
	\end{equation*}
	where $v=C_iC_jP_{i,j}(x,y)$, we choose $M_N(i+\lfloor N/2\rfloor)=C_i/4$. With this choice, we get $\hat{K}_{M_N}=\mathcal{K}$.

\appendix
	\section{Example for small distances}

	Let us consider $\bbd=1,N=3$ as the simplest example. We compute the polynomials that define the space of polynomials $V$ in Proposition \ref{propd1}.  We make $M_3(2),M_3(3)$ as free parameters and $M_3(1)=0$. We write  $K_{M_3}(d,d')$ for $1\le d\le d'\le 6$ in the following table.
	\begin{table}[H]
		\centering
		\begin{tabular}{|c|c|c|c|c|c|c|}
			\hline
			\diaghead{\hspace{3em}}{\vspace{-1mm}$d'$}{$d$}& $1$ & $2$ & $3$ & $4$ & $5$& $6$\\
			\hline
			$1$ & $\frac{1}{2}M_3(2)M_3(3)$ &  &  &  &  & \\
			\hline
			$2$ &  $-\frac{1}{2}M_3(2)M_3(3)$ & $0$ &  &  &  & \\
			\hline
			$3$ &  $-\frac{1}{2}M_3(2)M_3(3)$ & $0$ & $0$ &  &  & \\
			\hline
			$4$ & $0$ & $\frac{1}{2}M_3(2)^2$ & $0$ & $-\frac{1}{4}M_3(2)^2$ &  & \\
			\hline
			$5$ & $0$ & $\frac{1}{2}M_3(2)M_3(3)$ & $\frac{1}{2}M_3(2)M_3(3)$ & $0$ & $-\frac{1}{2}M_3(2)M_3(3)$ & \\
			\hline
			$6$ & $0$ & $0$ & $\frac{1}{2}M_3(3)^2$ & $0$ & $0$ & $-\frac{1}{4}M_3(3)^2$\\
			\hline
		\end{tabular}
		\caption{$K_{M_N}$ for $N=3$}
	\end{table}
	Then we compute $L_{M_3}(d_1,d_2)$ for $0\le d_1\le d_2\le 5$.
	\begin{table}[H]
		\centering
		\begin{tabular}{|c|c|c|}
			\hline
			\diaghead{\hspace{3em}}{\vspace{-1mm}$d_2$}{$d_1$}& $0$ & $1$\\
			\hline
			$0$ & $\frac{3}{4}\left( M_3(2)^2+M_3(3)^2 \right)$ &  \\
			\hline
			$1$ & $\frac{3}{4}\left( M_3(2)^2+M_3(3)^2\right)+\frac{1}{2}M_3(2)M_3(3)$ & $\frac{3}{4}\left( M_3(2)^2+M_3(3)^2 \right)+\frac{3}{2}M_3(2)M_3(3)$ \\
			\hline
			$2$ & $\frac{1}{4}M_3(2)^2+\frac{3}{4}M_3(3)^2+\frac{1}{2}M_3(2)M_3(3)$ & $\frac{1}{4}M_3(2)^2+\frac{3}{4} M_3(3)^2+M_3(2)M_3(3)$ \\
			\hline
			$3$ & $\frac{1}{4}\left( M_3(2)^2+M_3(3)^2 \right)+\frac{1}{2}M_3(2)M_3(3)$ & $\frac{1}{4}\left( M_3(2)^2+M_3(3)^2 \right)+\frac{1}{2}M_3(2)M_3(3)$ \\
			\hline
			$4$ & $\frac{1}{4}M_3(3)^2+\frac{1}{2}M_3(2)M_3(3)$ & $\frac{1}{4}M_3(3)^2 +\frac{1}{2}M_3(2)M_3(3)$ \\
			\hline
			$5$ & $\frac{1}{4}M_3(3)^2$ & $\frac{1}{4}M_3(3)^2$ \\
			\hline
		\end{tabular}
		\begin{tabular}{|c|c|c|}
			\hline
			\diaghead{\hspace{3em}}{\vspace{-1mm}$d_2$}{$d_1$}& $2$ & $3$\\
			\hline
			$2$ & $-\frac{1}{4}M_3(2)^2+\frac{3}{4}M_3(3)^2+\frac{1}{2}M_3(2)M_3(3)$ &\\
			\hline
			$3$ & $\frac{1}{4}\left( -M_3(2)^2+M_3(3)^2 \right)$ & $-\frac{1}{4}\left( M_3(2)^2+M_3(3)^2\right)-\frac{1}{2}M_3(2)M_3(3)$\\
			\hline
			$4$ & $\frac{1}{4}M_3(3)^2$ & $-\frac{1}{4}M_3(3)^2 -\frac{1}{2}M_3(2)M_3(3)$\\
			\hline
			$5$ & $\frac{1}{4}M_3(3)^2$ & $-\frac{1}{4}M_3(3)^2$\\
			\hline
		\end{tabular}
		\begin{tabular}{|c|c|c|}
			\hline
			\diaghead{\hspace{3em}}{\vspace{-1mm}$d_2$}{$d_1$}& $4$ & $5$\\
			\hline
			$4$ & $-\frac{1}{4} M_3(3)^2 -\frac{1}{2}M_3(2)M_3(3)$ &\\
			\hline
			$5$ & $-\frac{1}{4}M_3(3)^2$ & $-\frac{1}{4}M_3(3)^2$\\
			\hline
		\end{tabular}
		\caption{$L_{M_N}$ for $N=3$}
	\end{table}
	We define
	\begin{align*}
		P(x,y):=& \frac{3}{4}(U_2(x)+1)(U_2(y)+1)+\frac{1}{4}\left( U_4(x)+U_4(y)+U_4(x)U_2(y)+U_2(x)U_4(y)-U_4(x)U_4(y) \right)\\
		&+\frac{1}{4}\left( U_6(x)+U_6(y)+U_6(x)U_2(y)+U_2(x)U_6(y)-U_6(x)U_4(y)-U_4(x)U_6(y)-U_6(x)U_6(y) \right)\\
		Q(x,y):=& \frac{3}{2}U_2(x)U_2(y)+U_4(x)U_2(y)+U_2(x)U_4(y)\\
		&+\frac{1}{2}\left( U_2(x)+U_2(y)+U_4(x)+U_4(y)+U_6(x)+U_6(y)+U_8(x)+U_8(y) \right)\\
		&+\frac{1}{2}\left( U_6(x)U_2(y)+U_2(x)U_6(y)+U_8(x)U_2(y)+U_2(x)U_8(y)+U_4(x)U_4(y) \right)\\
		&-\frac{1}{2}\left( U_6(x)U_6(y)+U_8(x)U_6(y)+U_6(x)U_8(y)+U_8(x)U_8(y) \right)\\
		R(x,y):=& \frac{3}{4}\left( 1+U_2(x)+U_2(y)+U_2(x)U_2(y)+U_4(x)+U_4(y)+U_4(x)U_2(y)+U_2(x)U_4(y)+U_4(x)U_4(y) \right)\\
		&+\frac{1}{4}\left( U_6(x)+U_6(y)+U_6(x)U_2(y)+U_2(x)U_6(y)+U_6(x)U_4(y)+U_4(x)U_6(y)+U_6(x)U_6(y) \right)\\
		&+\frac{1}{4}\left( U_8(x)+U_8(y)+U_8(x)U_2(y)+U_2(x)U_8(y)+U_8(x)U_4(y)+U_4(x)U_8(y) \right)\\
		&-\frac{1}{4}\left( U_8(x)U_6(y)+U_6(x)U_8(y)+U_6(x)U_6(y)+U_8(x)U_8(y) \right)\\
		&+\frac{1}{4}U_{10}(x)\left( 1+U_2(y)+U_4(y)-U_6(y)-U_8(y)-U_{10}(y) \right)\\
		&+\frac{1}{4}U_{10}(y)\left( 1+U_2(x)+U_4(x)-U_6(x)-U_8(x)-U_{10}(x) \right)
	\end{align*}
We have the set of polynomials
\begin{equation*}
	V=\{C_1^2P(x,y)+C_1C_2Q(x,y)+C_2^2R(x,y)\mid C_1,C_2\in\R\}.
\end{equation*}
Proposition \ref{propd1} states that we can always find a perturbation to get the kernel in the form of
\begin{equation*}
	\sin^2(\pi k)\sin^2(\pi k')v(\cos(\pi k),\cos(\pi k'))\quad\text{for }v\in V.
\end{equation*}

\bibliographystyle{plain}

\bibliography{references}

\begin{thebibliography}{10}

\bibitem{basile2010energy}
G.~Basile, S.~Olla, and H.~Spohn.
\newblock Energy transport in stochastically perturbed lattice dynamics.
\newblock {\em Archive for Rational Mechanics and Analysis}, 195(1):171--203,
  2010.

\bibitem{borcea2016derivation}
L.~Borcea and J.~Garnier.
\newblock Derivation of a one-way radiative transfer equation in random media.
\newblock {\em Physical Review E}, 93(2):022115, 2016.

\bibitem{carrillo2022controlling}
J.~A. Carrillo, D.~Kalise, F.~Rossi, and E.~Tr{\'e}lat.
\newblock Controlling swarms toward flocks and mills.
\newblock {\em SIAM Journal on Control and Optimization}, 60(3):1863--1891,
  2022.

\bibitem{hannani2022wave}
A.~Hannani, M.~Rosenzweig, G.~Staffilani, and M.-B. Tran.
\newblock On the wave turbulence theory for a stochastic {KdV} type
  equation--generalization for the inhomogeneous kinetic limit.
\newblock {\em arXiv preprint arXiv:2210.17445}, 2022.

\bibitem{hannani2024controlling}
A.~Hannani, M.-B. Tran, M.-N. Phung, and E.~Tr\'{e}lat.
\newblock Controlling the rates of a chain of harmonic oscillators with a point
  langevin thermostat, 2024.

\bibitem{komorowski2023heat}
T.~Komorowski, J.L. Lebowitz, and S.~Olla.
\newblock Heat flow in a periodically forced, thermostatted chain.
\newblock {\em Communications in Mathematical Physics}, pages 1--45, 2023.

\bibitem{komorowski2022asymptotic}
T.~Komorowski and S.~Olla.
\newblock Asymptotic scattering by {Poissonian} thermostats.
\newblock In {\em Annales Henri Poincar{\'e}}, volume~23, pages 3753--3790.
  Springer, 2022.

\bibitem{KORH2020}
T.~Komorowski, S.~Olla, L.~Ryzhik, and H.~Spohn.
\newblock High frequency limit for a chain of harmonic oscillators with a point
  {Langevin} thermostat.
\newblock {\em Archive for Rational Mechanics and Analysis}, 237, 07 2020.

\bibitem{LS2007}
J.~Lukkarinen and H.~Spohn.
\newblock Kinetic limit for wave propagation in a random medium.
\newblock {\em Archive for Rational Mechanics and Analysis}, 183:93--162, 01
  2007.

\bibitem{piccoli2015control}
B.~Piccoli, F.~Rossi, and E.~Tr{\'e}lat.
\newblock Control to flocking of the kinetic {Cucker}--{Smale} model.
\newblock {\em SIAM Journal on Mathematical Analysis}, 47(6):4685--4719, 2015.

\bibitem{pomeau2019statistical}
Y.~Pomeau and M.-B. Tran.
\newblock {\em Statistical Physics of Non Equilibrium Quantum Phenomena},
  volume 967.
\newblock Springer Nature, 2019.

\bibitem{RPK96}
L.~Ryzhik, G.~Papanicolaou, and J.~Keller.
\newblock Transport equations for elastic and other waves in random media.
\newblock {\em Wave Motion}, 24:327--370, 12 1996.

\bibitem{staffilani2021wave}
G.~Staffilani and M.-B. Tran.
\newblock On the wave turbulence theory for a stochastic {KdV} type equation.
\newblock {\em arXiv preprint arXiv:2106.09819}, 2021.

\end{thebibliography}

\end{document}